\documentclass[12pt]{article}
\usepackage{epsfig,amssymb}

\makeatletter
\def\eqalign#1{\null\vcenter{\def\\{\cr}\openup\jot\m@th
  \ialign{\strut$\displaystyle{##}$\hfil&$\displaystyle{{}##}$\hfil
      \crcr#1\crcr}}\,}
\makeatother

\newcommand{\be}{\begin{equation}} 
\newcommand{\ee}{\end{equation}}
\newcommand{\beqa}{\begin{eqnarray}}
\newcommand{\eeqa}{\end{eqnarray}}
\newcommand{\bt}{\begin{theorem}}
\newcommand{\et}{\end{theorem}}
\newcommand{\bl}{\begin{lemma}}
\newcommand{\el}{\end{lemma}}
\newcommand{\bc}{\begin{corollary}}
\newcommand{\ec}{\end{corollary}}
\newcommand{\ba}{\begin{array}}
\newcommand{\ea}{\end{array}}

\newcommand{\la}{\label}
\newcommand{\ci}{\cite}

\newcommand{\wt}{\widetilde}
\newcommand{\de}{\delta}
\newcommand{\De}{\Delta}
\newcommand{\al}{\alpha}
\newcommand{\ga}{\gamma}
\newcommand{\Ga}{\Gamma}

\newcommand{\si}{\sigma}
\newcommand{\Si}{\Sigma}
\newcommand{\om}{\omega}

\newcommand{\lb}{\lambda}
\newcommand{\Lb}{\Lambda}
\newcommand{\ze}{\zeta}

\newcommand{\La}{\Lambda}
\newcommand{\ep}{\varepsilon }

\newcommand{\cD}{{\cal D}}

\newcommand{\bi}{\bibitem}

\def\bbr{\mathbb R}
\def\complex{\mathbb C}

\newtheorem{theorem}{Theorem}
\newtheorem{proposition}{Proposition}
\newtheorem{corollary}{Corollary}
\newtheorem{lemma}{Lemma}

\begin{document}
\begin{center}
{\Large\bf The Widom-Dyson constant for the gap probability in
random matrix theory}
\end{center}

%\maketitle

\vskip 1.2in
\centerline{ P.~ Deift}
\centerline{{\it Courant Institute of Mathematical Sciences,}}
\centerline{{\it New York, NY 10003, USA}}
\vskip .2in
\centerline{ A.~Its}
\centerline{{\it Department of Mathematical Sciences,}}
\centerline{{\it Indiana University -- Purdue University  Indianapolis}}
\centerline{{\it Indianapolis, IN 46202-3216, USA}}
\vskip .2in
\centerline{ I.~Krasovsky}
\centerline{{\it Department of Mathematical Sciences,}}
\centerline{{\it Indiana University -- Purdue University  Indianapolis}}
\centerline{{\it Indianapolis, IN 46202-3216, USA}}
\centerline{{\it and}}
\centerline{{\it Department of Mathematical Sciences}}
\centerline{{\it Brunel University}}
\centerline{{\it Uxbridge UB83PH}}
\centerline{{\it United Kingdom}}
\vskip .2in
\centerline{ X.~Zhou}
\centerline{{\it Department of Mathematics,}}
\centerline{{\it Duke University,}}
\centerline{{\it Durham, NC 27708-0320}}

\vskip .4in

%\centerline{{\it( Draft )}}

\newpage
\section{Introduction}
In this paper we consider an asymptotic question in the theory
of the Gaussian Unitary Ensemble of random matrices \ci{M}.
In the bulk scaling limit, the probability that there are no
eigenvalues in the interval $(0,2s)$ is given by $P_s=\det(I-K_s)$,
where $K_s$ is the trace-class operator with kernel
$$
K_s(x,y) =\frac{\sin(x-y)}{\pi(x-y)}
$$
acting on $L^2(0,2s)$. We are interested particularly in the behavior of
$P_s$ as $s\to\infty$.

In 1973, des Cloizeaux and Mehta \ci{dCM} showed that as $s\to\infty$
\be
\ln P_s=-{s^2\over 2} -{1\over 4}\ln s+ c +o(1),\la{as1}
\ee
for some constant $c$. In 1976, Dyson \ci{D} showed that $P_s$ in fact has a
full asymptotic expansion of the form
\be\la{dyson_as}
\ln P_s=-{s^2\over 2} -{1\over 4}\ln s+ c_0 +{a_1\over s}+ {a_2\over
s^2}+\cdots.
\ee
Dyson identified all the constants $c_0$, $a_1$, $a_2$, $\dots$. Of
particular interest is the constant $c_0$, which he found using earlier
work of Widom (see \ci{W1} and below) to be
\be
c_0={1\over 12}\ln 2+ 3\ze'(-1),\la{c0}
\ee
where $\ze(z)$ is the Riemann zeta-function.

The results in \ci{dCM} and \ci{D} were not fully rigorous. In \ci{W2},
Widom gave the first rigorous proof of the leading asymptotics in
(\ref{as1}) in the form
\be
\ln P_s=-{s^2\over 2}(1+o(1)).
\ee
In subsequent work \ci{W3,DIZ}, which also included the multi-interval
generalization, the form (\ref{dyson_as}) of the full asymptotic
expansion was verified rigorously, together with the correct constants
$a_1$, $a_2$, $\dots$. The expression (\ref{c0}) for the constant $c_0$,
however, remained unproven. This was because the methods in \ci{W2,W3} and
\ci{DIZ} naturally computed the asymptotics of $(d/ds)\ln P_s$, and the
constant of integration remained undetermined.

Recently, two proofs of (\ref{c0}) were given independently in the
literature in \ci{K} and \ci{E,BE}. The methods in the papers
\ci{K} and \ci{E,BE} are very different. Our goal in this paper is to
give a third proof of (\ref{c0}), which is closely related to the proof
in \ci{K}, but as explained below, does not rely on certain a priori
information. This means that our approach has the potential advantage of
being applicable to other problems involving the computation of critical
constants, where a priori information may not be available (see, e.g.,
\ci{DIK2}).

One way that one might try to evaluate $c_0$ is to express
\be\la{tr}
\ln P_s=\ln\det(I-K_s)=\int_0^1{d\over d\eta}{\rm tr\,}\ln(I-\eta K_s)d\eta
=-\int_0^1{\rm tr\,}((I-\eta K_s)^{-1}K_s)d\eta
\ee
and then evaluate ${\rm tr\,}((I-\eta K_s)^{-1}K_s)$ asymptotically as
$s\to\infty$ for each fixed $\eta\in(0,1)$ using steepest descent
methods as in \ci{Deift}, for example. However, it turns out that the
asymptotics of ${\rm tr\,}((I-\eta K_s)^{-1}K_s)$ as $s\to\infty$ have a
different form for $\eta<1$ and $\eta=1$. This means that one must
integrate the asymptotics in (\ref{tr}) over a boundary layer as
$\eta\to 1$, a difficult task which we have so far been unable to perform.
On the other hand, for $0<\ga<1$, we can indeed use (\ref{tr}) in the form
\[
\ln\det(I-\ga K_s)=-\int_0^\ga{\rm tr\,}((I-\eta K_s)^{-1}K_s)d\eta
\]
together with the Riemann-Hilbert/steepest-descent method to compute
the asymptotics of $\ln \det (I-\ga K_s)$ as $s\to\infty$, so 
reproducing the results in \ci{BT,BB}.

As mentioned above, Dyson's computation of $c_0$ in \ci{D} is based on 
an earlier calculation
of Widom \ci{W1}. In \ci{W1}, Widom considered, in particular, 
the Toeplitz determinant
$D_n(\al)$ with symbol given by the characteristic function of the
interval $(\al,2\pi-\al)$, $0<\al<\pi$. Thus $D_n(\al)=\det
(M_{i-j})_{i,j=0}^{n-1}$, where
$M_k=\int_\al^{2\pi-\al}e^{-ik\theta}d\theta/{2\pi}$, $k\in Z$. Widom
showed that for a fixed $\al$ as $n\to\infty$,
\be\la{widom_as}
\ln D_n(\al)=n^2\ln\cos{\al\over 2}
-{1\over 4}\ln\left(n\sin{\al\over 2}\right)+c_0+o(1),
\ee
where $c_0$ is the constant (\ref{c0}). What Dyson noted was that
for a fixed $s>0$,
\be\la{limD}
\lim_{n\to\infty} D_n\left({2s\over n}\right)=\det(I-K_s)=P_s
\ee
and hence, if the error term $o(1)$ in (\ref{widom_as}) was uniform as
$n\to\infty$, $\al\to 0$, $\al n\to\infty$, one could conclude from
(\ref{widom_as},\ref{limD}) that $c_0$ in (\ref{dyson_as}) is indeed
given by (\ref{c0}).
The main technical result in this paper, as in \ci{K}, is the proof
that the error term $o(1)$ is of the form $O(1/(n\sin{\al\over 2}))$,
which gives the desired uniformity.

Whereas $P_s$ is the gap probability for the Gaussian Unitary Ensemble
in the bulk scaling limit, we note that $D_n(\al)$ is the gap
probability for the Circular Unitary Ensemble \ci{M}. Formula
(\ref{limD}) is the scaling limit for this probability, 
and the fact that the limit also
gives $P_s$ is a well-known universality property.

In \ci{K}, the author uses steepest descent methods to show that for
$\ep>0$ fixed, there exists a (large) positive constant $s_0$ such that
\be\la{diff}
{d\over d\al}\ln D_n(\al)= -{n^2\over 2}\tan{\al\over 2}-{1\over 8}\cot
{\al\over 2}+O\left({1\over n\sin^2(\al/2)}\right)
\ee
for all $n>s_0$ and ${2s_0\over n}\le\al\le\pi-\ep$.
Integrating (\ref{diff}) over $(\al,\al_0)$, $2s_0/n\le\al<\al_0\le\pi-\ep$,
one obtains
\be\eqalign{
\ln D_n(\al)=\ln D_n(\al_0)+n^2\ln\cos{\al\over 2}-{1\over
4}\ln\sin{\al\over 2}-n^2\ln\cos{\al_0\over 2}+\\
{1\over 4}\ln\sin{\al_0\over 2}+O\left({1\over n\sin(\al/2)}\right).}\la{alal0}
\ee
Using Widom's result (\ref{widom_as}) for fixed $\al_0$, one obtains
for ${2s_0\over n}\le\al\le\pi-\ep$, $n>s_0$,
\be\la{asdelta}
\ln D_n(\al)=n^2\ln\cos{\al\over 2}-{1\over
4}\ln\left(n\sin{\al\over 2}\right)+c_0+
O\left({1\over n\sin(\al/2)}\right)+\de_n,
\ee
where $\de_n\to 0$ as $n\to\infty$. For any fixed $s>s_0$, one sets
$\al=2s/n$, and then using (\ref{limD}) and letting $n\to\infty$, one
obtains
\be
\ln P_s=-{s^2\over 2} -{1\over 4}\ln s+ c_0 +O\left({1\over s}\right),
\ee
which proves (\ref{c0}).

In this paper, we will derive an improved version of (\ref{alal0}),
viz.,
\be\la{asnodelta}
\ln D_n(\al)=n^2\ln\cos{\al\over 2}-{1\over
4}\ln\left(n\sin{\al\over 2}\right)+c_0+
O\left({1\over n\sin(\al/2)}\right),
\ee
for ${2s_0\over n}\le\al\le\pi-\ep$, $n>s_0$, where $s_0$ is again a
(large) positive constant.
Our proof of (\ref{asnodelta}) is direct and does not rely on Widom's
result (\ref{widom_as}). The proof is based on the following two principles:
\begin{enumerate}
\item[(i)] Asymptotics of $D_{n}(\alpha)$ as $\alpha \to \pi$ and $n$ is fixed.
\item[(ii)] Asymptotics of the solution of a regularized version of the \ci{DIZ}
Riemann-Hilbert problem (see below)
uniform for $2s_0/n \leq \alpha \leq \pi$.
\end{enumerate}

The solution of problem (i) is based in turn on the analysis of the
standard multiple-integral representation for $D_n(\al)$. 
The solution of problem (ii)
is based on a mapping of the original Riemann-Hilbert
problem posed on the arc $-\al\le\theta\le\al$ of the unit circle
to a problem on the fixed
interval $[-1, 1]$. The analysis then proceeds via the steepest descent
method for Riemann-Hilbert problems introduced by Deift and Zhou in \ci{DZ}
and further developed in \ci{DZpainleve2}, \ci{DVZzerodisp}, and also in
\ci{Dstrong}. This gives an asymptotic expression for the logarithmic derivative
$(d^2/d\al^2)\ln D_n(\al)$.
Formula (\ref{differential}) below together with its
integrated version (\ref{Dinteg}), plays a key role in this paper.

Note that in contrast to \ci{W2} and \ci{DIZ}, where the analysis of the
derivative $(d/ds)\ln P_s$ fails to identify the constant $c_0$, we may
now integrate $(d^2/d\al^2)\ln D_n(\al)$
from $\al\to\pi$, and the limit at $\pi$ is determined by
step (i). The result is the expression (\ref{asnodelta}).
By contrast, in \ci{W2} and \ci{DIZ} there is no convenient point $s_0$
from which we can integrate and then use to extract the relevant
asymptotics. The key device that makes our
method work is the $\Phi$-RH: In particular, we note that the 11-element in the
jump matrix for the $\Phi$-RH (see (\ref{Phidef}) et seq.) is {\it uniformly} small as
$n\to \infty$, for {\it all} $2s_0/n \leq \alpha \leq \pi$, $s_0 >> 1$, and for
all $\lambda$ in a compact subset of $(-1,1)$. It is this uniformity in $\alpha$
as $n \to \infty$ that makes it possible to control the integration from
$\alpha = \pi$ to $\alpha=2s_0/n$.

In Section 2 we analyze step (i), and in section 3, step (ii). Finally,
in Section 4, we prove (\ref{asnodelta}).

\section{Step (i). Multiple integral analysis.}
For the analysis of $D_{n}(\alpha)$ as $\alpha \to \pi$
we use the multiple integral (e.g., \ci{Szego, Deift})
\begin{equation}\label{multint}
D_{n}(\alpha) = \frac{1}{(2\pi)^{n}n!}
\int_{C_\alpha}\cdots\int_{C_\alpha}
\prod_{1\leq j < k \leq n}|e^{i\theta_{j}} - e^{i\theta_{k}}|^{2}
d\theta_{1}\dots d\theta_{n},
\end{equation}
where $C_\al$ is the arc $\al\le\theta\le 2\pi-\al$ of the unit circle.
The integrals are taken from $\al$ to $2\pi-\al$.
Setting
$$
\alpha = \pi -\beta, \quad \beta > 0,
$$
and
$$
\theta_{j} = \pi + \beta x_{j},
$$
we rewrite (\ref{multint}) as follows:
\begin{equation}\label{multint2}
D_{n}(\alpha) = \frac{1}{(2\pi)^{n}n!}\beta^{n}
\int_{-1}^{1}\cdots\int_{-1}^{1}
\prod_{1\leq j < k \leq n}|e^{i\beta x_{j}} - e^{i\beta x_{k}}|^{2}
dx_{1}\dots dx_{n}.
\end{equation}
Observe that
$$
\prod_{1\leq j < k \leq n}|e^{i\beta x_{j}} - e^{i\beta x_{k}}|^{2}
= \beta^{n(n-1)}
\left(\prod_{1\leq j < k \leq n}| x_{j} -  x_{k}|^{2} + O_n(\beta^2)\right).
$$
Hence we arrive at the relation
\begin{equation}\label{multint3}
D_{n}(\alpha) = \frac{1}{(2\pi)^{n}n!}\beta^{n^2}
\left \{\int_{-1}^{1}\cdots\int_{-1}^{1}
\prod_{1\leq j < k \leq n}|x_{j} - x_{k}|^{2}
dx_{1}\dots dx_{n} + O_n(\beta^2)\right\}.
\end{equation}
The multiple integral in this formula can be expressed in terms
of the norms of the Legendre polynomials. Indeed (see, e.g., \ci{Szego})
\be
\int_{-1}^{1}\cdots\int_{-1}^{1}
\prod_{1\leq j < k \leq n}|x_{j} - x_{k}|^{2}
dx_{1}\dots dx_{n}
= n!\prod_{k=0}^{n-1}h_{k},
\ee
where $h_{n}$ are the normalization constants of the
monic polynomials orthogonal on the interval $[-1, 1]$
with the unit weight:
$$
p_{n}(x) = x^n +\dots, \quad \int_{-1}^{1}p_{n}(x)p_{m}(x)dx
=h_{n}\de_{nm}.
$$
Let $P_{n}(x)$ denote the standard Legendre polynomials \ci{Szego}. Since
$$
P_{n}(x) = \frac{(2n)!}{2^n(n!)^2}x^n +\dots
\quad \mbox{and}\quad \int_{-1}^{1}P^2_{n}(x)dx = \frac{2}{2n+1}
$$
we conclude that
$$
p_{n}(x) = \frac{2^n(n!)^2}{(2n)!}P_{n}(x)
$$
and
$$
h_{n} =  \int_{-1}^{1}p^2_{n}(x)dx =  \frac{2^{2n}(n!)^4}{[(2n)!]^2}
\int_{-1}^{1}P^2_{n}(x)dx = \frac{2^{2n}(n!)^4}{[(2n)!]^2}\frac{2}{2n+1}.
$$
This leads us to the following representation of $D_{n}(\alpha)$ in the
neighborhood of $\alpha =\pi$:
\begin{equation}\label{Dn1}
\ln D_{n}(\alpha) = n^2\ln \beta - n\ln 2\pi +\ln A_{n} + O_n(\beta^2),
\quad \alpha = \pi - \beta,
\end{equation}
where
\begin{equation}\label{An}
A_{n} = \prod_{k=0}^{n-1}h_{k} = \prod_{k=0}^{n-1}
\frac{2^{2k}(k!)^4}{[(2k)!]^2}\frac{2}{2k+1}.
\end{equation}
For later reference, note that the asymptotic relation (\ref{Dn1}) is clearly
differentiable, for fixed $n$, with respect to $\alpha$.
Also, for fixed $n$, the term $O_n(\beta^2)\to 0$ as $\beta=\pi-\al\to 0$;
no claim is made here about the behavior of $O_n(\beta^2)$ as $n\to\infty$.

Widom's constant,
$c_{0} = \frac{1}{12}\ln 2 + 3\zeta'(-1)$, is generated by the
quantity $A_{n}$. In fact, it is shown in \ci{W1}, using results from 
classical analysis, that
\begin{equation}\label{Anas}
A_{n} = e^{c_{0}}n^{-1/4}(2\pi)^n2^{-n^2}(1 + o(1)),\quad n \to \infty.
\end{equation}
The appearance of the zeta function
is due to the presence of the products of factorials. Indeed,
$$
\ln \prod_{k=1}^{n}k! = \sum _{k=1}^{n}\ln k!
=(n+1)\ln n! - \sum_{k=1}^{n}k\ln k,
$$
and one can expect that the asymptotics of the sum on the r.h.s.
of the last equation is
related to $\zeta'(-1)$. The exact relation (see again \ci{W1})
reads as follows:
$$
\sum_{k=1}^{n}k\ln k = \left(\frac{1}{2}n^2 +\frac{1}{2}n +
\frac{1}{12}\right)
\ln n - \frac{1}{4}n^2 +\frac{1}{12} - \zeta'(-1) + o(1).
$$
Applying this formula and the asymptotics of the Gamma-function to (\ref{An}) 
yields (\ref{Anas}).

\section{Step 2. Riemann-Hilbert analysis.}
Denote the complement of $C_\al$ in the unit circle by
$\Gamma_{\alpha}=\{-\alpha < \theta < \alpha\}$
traversed counterclockwise (see Figure 1).
Let  $m(z) \equiv
m(z;n,\alpha)$ be the solution of the following $2\times 2$ Riemann-Hilbert
problem  posed on $\Gamma_{\alpha}$:
\begin{itemize}
\item $m(z)$ is holomorphic for all $z\notin \overline{\Gamma_{\alpha}}$
\item $m(\infty) = I$
\item $ m_{-}(z) = m_{+}(z)\pmatrix{2 & -z^n \cr\cr
      z^{-n} & 0}, \quad z \in \Gamma_{\alpha}$
\end{itemize}
Here, as usual, $m_+(z)$ (respectively, $m_-(z)$) are the $L^2$ boundary values of 
$m(z')$ as $z'\to z\in\Ga_\al$ non-tangentially from the ``$+$'' side $\{|z|<1\}$ 
(respectively, ``$-$'' side $\{|z|>1\}$).
We shall refer to this Riemann-Hilbert problem as the ``$m$-RH problem''.

\begin{figure}
\centerline{\psfig{file=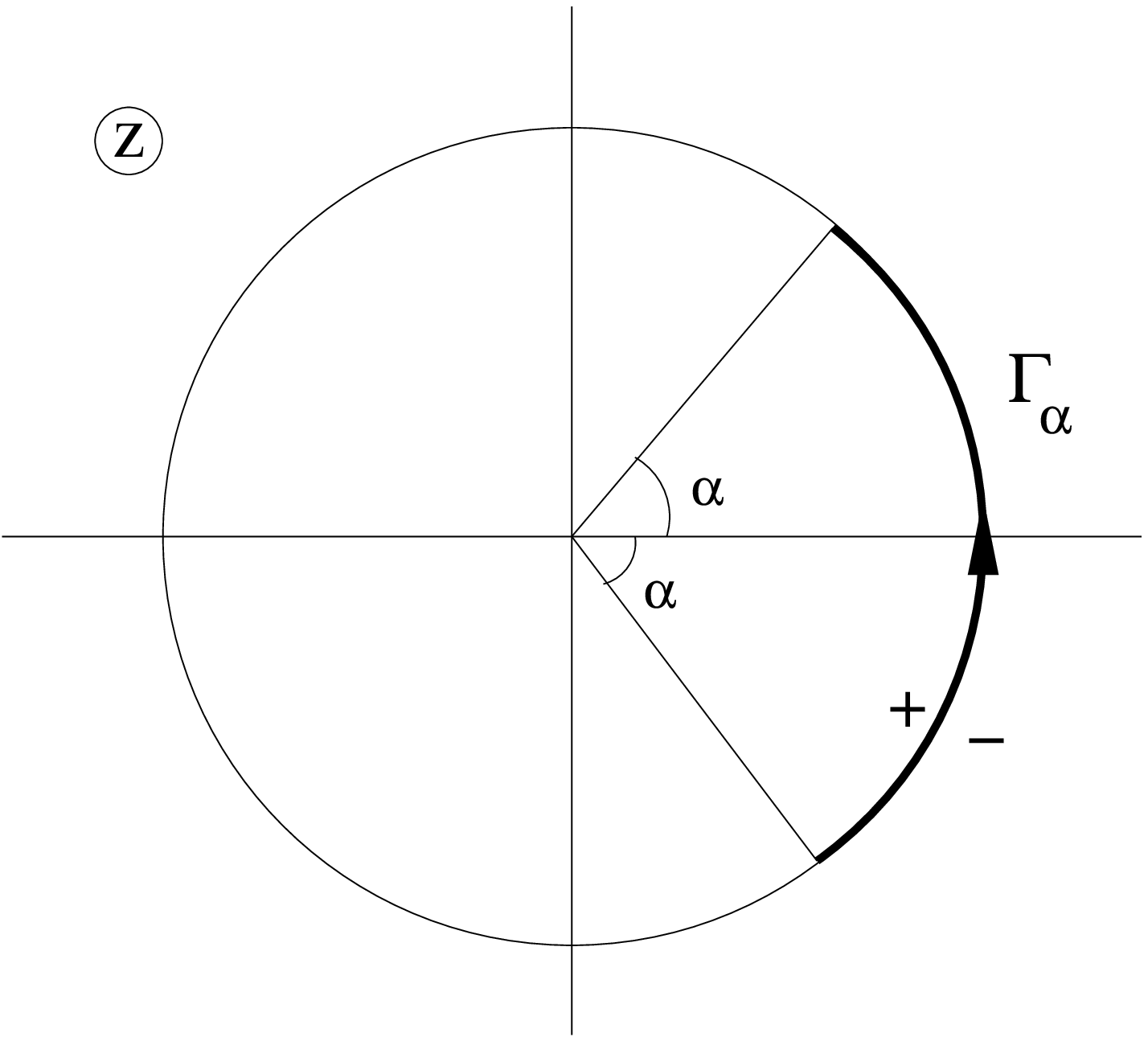,width=3.0in,angle=0}}
\vspace{0cm}
\caption{
Contour for the $m$-RH problem.}
\label{fig1}
\end{figure}

{\theorem\label{theor1} Let $0 < \alpha < \pi$ and $n >0$.
Then the $m$-RH problem has a (unique) solution $m(z; n, \alpha)$, 
and the Toeplitz determinant $D_{n}(\alpha)$ is related to 
$m(z; n, \alpha)$ by the following differential and difference identities:
\begin{equation}\label{difference}
\frac{D_{n+1}(\alpha)}{D_{n}(\alpha)}
= m_{11}(0;n, \alpha),
\end{equation}
\begin{equation}\label{differential}
\frac{d^2}{d\alpha^2}\ln D_{n}(\alpha)
= -\frac{n^2}{\sin^{2}\alpha}[m_{12}(0;n, \alpha)]^2.
\end{equation}}

\noindent
{\bf Remark} 
As a function of $z$, $m(z; n, \alpha)$ has a continuous extension up
to the boundary $\overline{\Gamma_{\alpha}}$, apart from the two end points
$e^{i\al}$ and $e^{-i\al}$, where it has logarithmic singularities. Moreover,
$m_{\pm}(z)$ admit analytic continuations into a neighborhood of
every point $z$ of the open arc $\Ga_\al=\overline{\Ga_\al}\setminus\{e^{i\al},
e^{-i\al}\}$.
Note also that $\det m(z; n, \alpha)=1$ by a standard calculation.
These properties of $m(z; n, \alpha)$ are inherited by solutions of 
the transformed Riemann-Hilbert problems introduced below.

Theorem 1  was proved in \ci{DIZ}\footnote{
There are some differences from the notation in \ci{DIZ}, namely,
our contour is $\Gamma_\alpha$ instead of $C$ ($\Ga_\al$ rotated by
$\pi$) in \ci{DIZ},
and we make the following choice for the functions $f_i$, $g_i$
which build up the kernel: $f_1=z^{n/2}$, $f_2=z^{-n/2}$,
$g_1=z^{-n/2}/(2\pi i)$, $g_2=-z^{n/2}/(2\pi i)$.}
(cf. Eqs (6.14) and (6.82))
using standard techniques from the
theory of integrable systems: derivation of the relevant Lax pair,
identification of $D_{n}(\alpha)$ as the relevant tau-function
etc. The differential identity (\ref{differential})
will be of central importance for the analysis below.

A standard calculation shows that the $m$-RH problem has no solution for
$\al=\pi$. However, as we now demonstrate, 
the $m$-RH problem can be
regularized for all $\al$ in the range, including $\al=\pi$, 
by a simple sequence of transformations.

\subsection{Mapping onto a fixed interval.}
For $0<\al<\pi$, the linear-fractional transformation,
\begin{equation}\label{z_lambda}
\lambda = -i\cot{\frac{\alpha}{2}}\, \,\frac{z-1}{z+1},
\quad z = \frac{1 + i\lambda \tan{\frac{\alpha}{2}}}
{1 - i\lambda \tan{\frac{\alpha}{2}}},
\end{equation}
maps the arc $\Gamma_{\alpha}$ onto the interval
$(-1,1)$ and transforms the $m$-RH problem to the following
Riemann-Hilbert problem posed on the interval $(-1,1)$
traversed from $-1$ to $1$ (see Figure 2):
\begin{itemize}
\item $Y(\lambda)$ is holomorphic for all $\lambda\notin [-1, 1]$
\item $Y(\infty) = I$
\item $ Y_{-}(\lambda) = Y_{+}(\lambda)\pmatrix{2 & -
\left(\frac{1 + i\lambda \tan{\frac{\alpha}{2}}}
{1 - i\lambda \tan{\frac{\alpha}{2}}}\right)^n \cr
 \left(\frac{1 + i\lambda \tan{\frac{\alpha}{2}}}
{1 - i\lambda \tan{\frac{\alpha}{2}}}\right)^{-n}& 0},
\quad \lambda \in (-1, 1)$
\end{itemize}
We shall refer to this Riemann-Hilbert problem as the
``$Y$-RH problem''. The relation between the $Y$-RH problem
and  the original $m$-RH problem is given by the equation
\begin{equation}\label{Y_m}
m(z;n, \alpha) = Y^{-1}(-i\cot{\frac{\alpha}{2}};n, \alpha)
Y(\lambda(z);n,\alpha).
\end{equation}

\begin{figure}
\centerline{\psfig{file=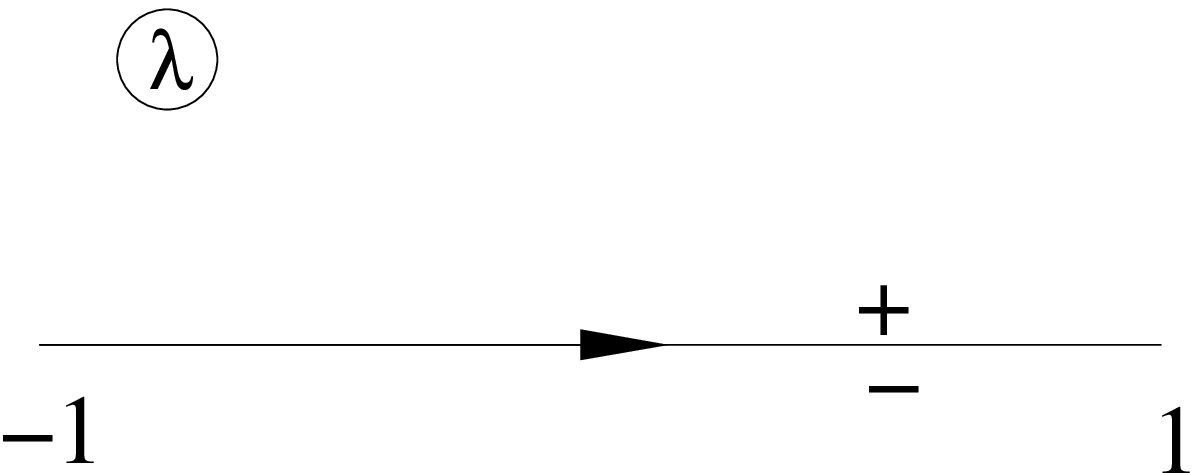,width=3.0in,angle=0}}
\vspace{0cm}
\caption{
Contour for the $Y$-RH problem.}
\label{fig2}
\end{figure}

The  $Y$-RH problem is still irregular at $\alpha = \pi$. Indeed,
the function 
$\phi(\lambda, \alpha)\equiv\frac{1 + i\lambda \tan(\alpha/2)}
{1 - i\lambda \tan(\alpha/2)}$ is discontinuous at
$(\lambda, \alpha) = (0,\pi)$.
We have for the jump matrix $J_{Y}(\lambda, \alpha)$ of $Y(\lb)$ as $\al\to\pi$:
$$
J_{Y}(\lambda, \pi) = \cases{
\pmatrix{2& -(-1)^{n}\cr (-1)^{n}&0}
,& $\lb \neq 0$\cr
\pmatrix{2& -1\cr 1 &0}
,& $\lb = 0$,}
$$
which demonstrates the difficulty 
for odd $n$. For even $n$, however, the jump matrix $J_{Y}(\lambda, \pi)$ is
continuous and constant throughout the whole interval $(-1, 1)$.
This implies the solvability of the $Y$-RH problem
at $\alpha = \pi$ for even $n$; in fact,
one easily checks that
\[
Y(\lb;n=2k,\pi)=\pmatrix{1&1\cr 1&0}
\pmatrix{1&{-1\over 2\pi i}\int_{-1}^1{ds\over s-\lb}\cr 0&1}
\pmatrix{1&1\cr 1&0}^{-1}.
\]
However, regardless of the parity of $n$,
the convergence of $J_{Y}(\lambda, \alpha)$ to $J_{Y}(\lambda, \pi)$
is not uniform in $\lb$, and this creates a significant difficulty
in the direct analysis of the behavior of the solution $Y(\lb;n,\alpha)$ 
near $\alpha = \pi$.
   
As we now show (see the $\Phi$-RH problem below), 
the  $Y$-problem can be regularized by performing one more
step which is familiar in the formalism of the nonlinear steepest descent
method.

\subsection{$g$-function transformation.}

Following the nonlinear steepest descent method for Riemann-Hilbert
problems (see, e.g., \ci{Dstrong}), we introduce the following ``g''
function:
\begin{equation}\label{g_def}
g(\lambda) \equiv \frac{1 + i\sqrt{\lambda^2 - 1}\sin{\frac{\alpha}{2}}}
{1 + i\lambda \tan{\frac{\alpha}{2}}}.
\end{equation}
This is essentially the $g$-function of section 6 of \ci{DIZ}
written in the variable $\lambda$ (see equation (\ref{gz})
below). It possesses the
following characteristic properties:
\begin{description}
\item [{\bf(a)}] $g(\lambda)$ is holomorphic for all $\lambda\notin [-1, 1]$.
Here we fix the square root by the condition
$$
\sqrt{\lambda^2 - 1} \sim \lambda, \quad \lambda \to \infty.
$$
\item [{\bf(b)}] $g(\lambda) \neq 0$ for all $\lambda\notin [-1, 1]$.
At the points $\lambda = -i\cot{\frac{\alpha}{2}}$
(or $z = \infty$) and  $\lambda = i\cot{\frac{\alpha}{2}}$
(or $z = 0$) the values of the function $g(\lambda)$ are:
\begin{equation}\label{z_infty_0}
g(-i\cot{\frac{\alpha}{2}}) = 1\quad \mbox{and} \quad
g(i\cot{\frac{\alpha}{2}}) = \cos^{2}{\frac{\alpha}{2}} \equiv \kappa.
\end{equation}
\item [{\bf(c)}] The boundary values of $g_{\pm}(\lambda)$,
$\lambda \in [-1, 1]$ satisfy the following equations:
\begin{equation}\label{g_bound1}
g_{+}g_{-} = \kappa \frac{1 - i\lambda \tan{\frac{\alpha}{2}}}
{1 + i\lambda \tan{\frac{\alpha}{2}}}
\end{equation}
and
\begin{equation}\label{g_bound2}
\frac{g_{+}}{g_{-}} =  \frac{1 - \sqrt{1-\lambda^2}
\sin{\frac{\alpha}{2}}}
{1 + \sqrt{1 -\lambda^2} \sin{\frac{\alpha}{2}}}.
\end{equation}
\item [{\bf(d)}] The behavior of $g(\lambda)$ as $\lambda \to \infty$
is described by the asymptotic relation
\begin{equation}\label{ginfty}
g(\lambda) = \cos{\frac{\alpha}{2}} + O\left(\frac{1}{\lambda}\right).
\end{equation}
\end{description}

It is worth noticing that
\begin{equation}\label{gz}
g(\lambda(z)) = \frac{z+1+\sqrt{(z-e^{i\alpha})(z-e^{-i\alpha})}}
{2z} \equiv \varphi(z).
\end{equation}
If one changes $1$ to $-1$ in the numerator, then $\varphi(z)$
becomes the $g$-function of section 6 of \ci{DIZ}. The change of sign
is due to the fact that the Riemann-Hilbert problem considered
in \ci{DIZ} is defined on the arc $C=e^{i\pi}\Ga_\al$
rather than on $\Ga_\al$ (cf. footnote 1 above).

Equation (\ref{g_bound2}) has an important consequence.
Fix $0<\de<1$ and $0<\al\le\pi$. Then
the following inequality holds:
\begin{equation}\label{ineq}
\left |\frac{g_{+}}{g_{-}}\right | \leq \ep_{0} < 1,\qquad
\lambda \in [-1 + \de, 1 -\de],
\end{equation}
for some $\ep_0=\ep_0(\de,\al)>0$. Of course, for all $\lb\in (-1,1)$
and $\al\in(0,\pi)$, we have
\[
\left |\frac{g_{+}}{g_{-}}\right | \leq 1.
\]

Following the steepest descent method, we transform the original
Riemann-Hilbert problem by the formula
\begin{equation}\label{Phidef}
Y(\lambda) \mapsto \Phi(\lambda) \equiv Y(\lambda)g^{-n\sigma_{3}}
\kappa^{\frac{n}{2}\sigma_{3}},
\end{equation}
where $\sigma_{3} = \pmatrix{1 & 0\cr 0 &-1}$ is the third Pauli matrix. From the 
properties of the $g$-function listed above, it follows that
the matrix function $\Phi(\lambda) \equiv \Phi(\lambda;n, \alpha)$
is the solution of the following Riemann-Hilbert problem,
which we shall refer to as the ``$\Phi$-RH problem'':
\begin{itemize}
\item $\Phi(\lambda)$ is holomorphic for all $\lambda\notin [-1, 1]$
\item $\Phi(\infty) = I$
\item $ \Phi_{-}(\lambda) = \Phi_{+}(\lambda)
\pmatrix{2 \left[\frac{1 - \sqrt{1-\lambda^2}
\sin{\frac{\alpha}{2}}}
{1 + \sqrt{1 -\lambda^2} \sin{\frac{\alpha}{2}}}\right]^{n}
& -1 \cr\cr
1 & 0}, \quad \lambda \in (-1, 1)$
\end{itemize}

In view of (\ref{Y_m}), the original function $m(z)$ is
related to the solution $\Phi(\lambda)$ by the formulae:
$$
m(z;n,\alpha) = \kappa^{\frac{n}{2}\sigma_{3}}
\Phi^{-1}(-i\cot{\frac{\alpha}{2}};n,\alpha)
\Phi(\lambda(z);n,\alpha)g^{n\sigma_{3}}(\lambda(z))
\kappa^{-\frac{n}{2}\sigma_{3}}
$$
\begin{equation}\label{Phim}
= \kappa^{\frac{n}{2}\sigma_{3}}
\Phi^{-1}(-i\cot{\frac{\alpha}{2}};n,\alpha)
\Phi(\lambda(z);n,\alpha)\varphi^{n\sigma_{3}}(z)
\kappa^{-\frac{n}{2}\sigma_{3}}.
\end{equation}

As indicated earlier, the $\Phi$-RH problem is regularized. Indeed,
note first that the jump matrix for the $\Phi$-RH problem is now continuous for all 
$\lb\in[-1,1]$ and $\al\in [0,\pi]$ with the end point $\al=\pi$ included.
Moreover, for all $0\le\al\le\pi$, the $\Phi$-RH problem is (uniquely) $L^2$-solvable,
by the following argument. Consider the $\Phi$-RH problem as defined on the whole real line
with discontinuities at $\pm 1$ (the jump matrix outside $(-1,1)$ is the identity). 
The limiting value of the jump matrix at these two points 
from inside the interval $(-1,1)$ is $\pmatrix{ 2&-1\cr 1&0}$, whose only eigenvalue is 
$1\notin (-\infty,0]$. By Theorem 5.16 of \ci{LS}, the RH problem is $L^2$ Fredholm.
Now by (the proof of) Theorem 9.3 in \ci{Zhou}, 
a Fredholm RH problem with a jump matrix $v$ on $\bbr$
is $L^2$-solvable if $v+v^*\ge 0$ everywhere, and $v+v^*>0$ on 
a set of positive Lebesgue measure.
These conditions are clearly satisfied in our case. Therefore the $\Phi$-RH problem is
$L^2$-solvable.

Theorem \ref{theor1} and equation (\ref{Phim}) yield
representations of the Toeplitz determinant $D_{n}(\alpha)$
in terms of the solution of the $\Phi$-RH problem:
\begin{eqnarray}\label{differencePhi}
\frac{D_{n+1}(\alpha)}{D_{n}(\alpha)}
= \Theta(n,\alpha)\cos^{2n}{\frac{\alpha}{2}},\\
\label{differentialPhi}
\frac{d^2}{d\alpha^2}\ln D_{n}(\alpha)
= -\frac{n^2}{\sin^{2}\alpha}\Delta(n,\alpha),
\end{eqnarray}
where
$$
\Theta(n,\alpha) = \left[\Phi^{-1}(-i\cot{\frac{\alpha}{2}};n,\alpha)
\Phi(i\cot{\frac{\alpha}{2}};n,\alpha)\right]_{11}
$$
and
\begin{equation}\label{Deltadef}
\Delta(n,\alpha) = \left[\Phi^{-1}(-i\cot{\frac{\alpha}{2}};n,\alpha)
\Phi(i\cot{\frac{\alpha}{2}};n,\alpha)\right]^{2}_{12}.
\end{equation}

\subsection{Asymptotic analysis of the $\Phi$-RH problem.}

By standard arguments, using inequality
(\ref{ineq}), one expects that $\Phi(\lambda)$
is approximated by the function
\begin{equation}\label{Phiinfty}
N(\lambda) =
\pmatrix{\frac{\beta(\lambda) + \beta^{-1}(\lambda)}{2} &
\frac{\beta(\lambda) - \beta^{-1}(\lambda)}{2i}\cr\cr
-\frac{\beta(\lambda) - \beta^{-1}(\lambda)}{2i} &
\frac{\beta(\lambda) + \beta^{-1}(\lambda)}{2}},
\quad \beta(\lambda) = \left(\frac{\lambda -1}
{\lambda +1}\right)^{\frac{1}{4}},
\end{equation}
$$
\beta(\infty) = 1,
$$
which solves the
model Riemann-Hilbert problem:
\begin{itemize}
\item $N(\lambda)$ is holomorphic for all $\lambda\notin [-1, 1]$
\item $N(\infty) = I$
\item $N_{-}(\lambda) = N_{+}(\lambda)
\pmatrix{0 & -1 \cr\cr
1 & 0}, \quad \lambda \in (-1, 1)$
\end{itemize}
In order to estimate the precision of this approximation, we need to
consider the $\Phi$-RH problem for $\lb$ near $\pm 1$.
The following result, which allows for complex values of $\al$ in 
a neighborhood of $\al=\pi$, is basic for our analysis 
(the need for this complex extension
will be apparent towards the end of the paper, see (\ref{cauchy}) below).

{\theorem \label{theorem2} Let $\de$
be a positive number less than $1/4$.
Introduce the domain
$$
\Omega^{(\de)} =\complex\setminus
(\overline{U\cup \wt U}),
$$
were $U$ ($\wt U$) denotes the open disk of radius $\de$ 
centered at $1$ (respectively, $-1$). Let also $\ep$
be a positive number less than $\pi - 2$ and denote $\cD_\ep(\pi)$ 
the disk in the $\al$-plane of radius 
$\ep$ centered at $\al=\pi$.
Set
$$
\rho=n\left|\sin{\al\over 2}\right|.
$$
Then, for $\de$ and $\ep$ sufficiently small, 
there exists $s_{0} > 0$ such that for all
$ \al\in\left[\frac{2s_0}{n},\pi-\ep\right]\cup\cD_\ep(\pi),$
and $n \geq s_{0}$, the solution of the $\Phi$-RH
problem exists (and is unique) and  satisfies the estimate
\begin{equation}\label{th21}
\Phi(\lambda;n,\alpha) = \left(I + O\left(\frac{1}{(1
+|\lambda|)\rho}\right)\right) N(\lambda), \quad \rho \to \infty,
\end{equation}
uniformly for $\lambda \in \Omega^{(2\de)}$ and 
$ \al\in\left[\frac{2s_0}{n},\pi-\ep\right]\cup\cD_\ep(\pi).$
Moreover, this estimate can be extended to a full
asymptotic series in inverse powers of $\rho$;
in particular, the order $\rho^{-3}$ extension of (\ref{th21})
reads:
\begin{equation}\label{th2105}
\Phi(\lambda;n,\alpha) = \left(I + 
R_1(\lb) + R_{2}(\lb) + R_r(\lb)\right)
N(\lambda), 
\end{equation}
where
\begin{eqnarray}
R_1(\lb)={1\over 16i\rho}\left[
{1\over 1-\lb}\pmatrix{-1&i\cr i&1}+
{1\over 1+\lb}\pmatrix{1&i\cr i&-1}\right],\qquad
\lb\in\complex\setminus(\overline{U\cup\wt U}),\la{R105}\\
R_2(\lb)=
{1\over 2^8\rho^2}\left[
{1\over 1-\lb}\pmatrix{1&8i\cr -8i&1}+
{1\over 1+\lb}\pmatrix{1&-8i\cr 8i&1}\right],\qquad
\lb\in\complex\setminus(\overline{U\cup\wt U}),\la{R205}
\end{eqnarray}
\be
R_r(\lb)=
O\left(\frac{1}{(1 +|\lambda|)\rho^{3}}\right),
\quad \rho \to \infty,\la{Rr05}
\ee
uniformly for $\lambda \in \Omega^{(2\de)}$ and 
$\al\in\left[\frac{2s_0}{n},\pi-\ep\right]\cup\cD_\ep(\pi)$.
}

{\bf Remark 1} The last statement (\ref{Rr05}) 
means that there exist positive constants $C$ and $s_{0}$,
depending on $\ep$ and $\de$ only,  such that
\begin{equation}\label{meaning05}
|R_{r}(\lambda)| \leq \frac{C}{(1 +|\lambda|)\rho^{3}},
\end{equation}
$$
\forall \lambda \in \Omega^{(2\de)}, \quad 
\forall \al\in\left[\frac{2s_0}{n},\pi-\ep\right]\cup\cD_\ep(\pi),
\quad  \forall n: s_0\leq n.
$$
We shall also assume that $\ep$ is small enough for
the inequality, 
\begin{equation}\label{ineqsin}
\left|\sin {\al\over 2}\right| \geq \frac{1}{2}
\end{equation}
to take place
for all $\alpha \in \cD_\ep(\pi)$, and hence 
\begin{equation}\label{rhos05}
\rho \geq \frac{1}{2} s_{0},
\end{equation}
for all $ \al\in\left[\frac{2s_0}{n},\pi-\ep\right]\cup\cD_\ep(\pi)$
and $s_{0} \leq  n$. 

{\bf Remark 2} 
Part of the assertion of Theorem 2 is that the solution of the $\Phi$-RH problem 
exists and is unique for all $\alpha\in [2s_0/n,\pi-\ep ]\cup\cD_\ep(\pi)$ and $n\ge s_0$
with $s_0$ sufficiently large. This is all we need in the analysis that follows;
however, the solution of the $\Phi$-RH problem actually exists and is unique for all
$\al\in[0,\pi-\ep]\cup\cD_\ep(\pi)$ and {\it all} $n>0$ for some (possibly smaller)
$\ep>0$. Indeed, by the discussion following (\ref{Phim}),
the $\Phi$-RH problem is solvable for all $\al\in[0,\pi]$, $n>0$, and also 
for all $\al\in\cD_{\ep'}(\pi)$, $0<n<s_0$ for some $\ep'>0$ by continuity of the 
jump matrix at $\al=\pi$.
By Theorem 2, the $\Phi$-RH problem is solvable for all $\al\in\cD_{\ep}(\pi)$,
$n\ge s_0$. 
Thus the $\Phi$-RH problem is solvable for all
$n>0$ on $[0,\pi-\ep_1]\cup\cD_{\ep_1}(\pi)$, where $\ep_1=\min(\ep,\ep')$.

%Although we will not use this fact,
%the $\Phi$-RH problem is actually solvable
%for all $\alpha\in (0,\pi-\ep ]\cup\cD_\ep(\pi)$ and $n >0$
%(and sufficiently small $\ep$).
%This is a consequence of the solvability of
%the original $m$-RH problem: see Theorem \ref{theor1}.
%In more detail, the solvability
%of the $m$-RH problem for $0 < \alpha < \pi$ is equivalent
%to the inequality $D_{n}(\alpha) \neq 0$ (cf. \cite{DIZ}) 
%which follows, e.g., from the integral representation (\ref{multint2}).
%For $\alpha \in  \cD_\ep(\pi)$ we shall use a modification
%of the formula (\ref{multint2}), namely
%\begin{equation}\label{multint2mod}
%D_{n}(\alpha) = \frac{1}{(2\pi)^{n}n!}\beta^{n}
%\int_{-1}^{1}\cdots\int_{-1}^{1}
%\prod_{1\leq j < k \leq n}(e^{i\beta x_{j}} - e^{i\beta x_{k}})
%(e^{-i\beta x_{j}} - e^{-i\beta x_{k}})
%dx_{1}\dots dx_{n}.
%\end{equation}
%This equation provides an analytic continuation of $D_{n}(\alpha)$
%in the disc $\cD_\ep(\pi)$, and shows that $D_{n}(\alpha) \neq 0$
%for $\alpha \in \cD_\ep(\pi)\setminus \{\pi\}$. This, by virtue of
%the relations (\ref{multint3}) and (\ref{differentialPhi}), yields
%the solvability of the $\Phi$-RH problem for all 
%$\alpha\in (0,\pi-\ep ]\cup\cD_\ep(\pi)$ and $n >0$. It also worth noticing
%that the orthogonal polynomial representation of the solution of $m(z;n, \alpha)$ 
%of the $m$-RH problem (see \cite{deift}, \cite{bdj}) implies a similar
%representation of the solution of the $\Phi$-RH problem.  

{\bf Remark 3} The local analyticity of the jump matrix of the $\Phi$-RH problem
implies that both boundary values of the
function $\Phi(\lambda)$ on $(-1,1)$, i.e. the functions $\Phi_{\pm}(x)$,
admit the analytic  continuation in the neighborhood of every point
of the interval $(-1 + \de, 1 - \de)$.

{\bf Proof of Theorem \ref{theorem2}.}
We shall now construct parametrices in $U$ and $\wt U$ which are
solutions of the $\Phi$-RH problem in these neighborhoods
with the condition at infinity replaced by the requirement that they
match $N(\lb)$ at the disks' boundaries to leading order
(cf. \ci{DIZ,K}).

\begin{figure}
\centerline{\psfig{file=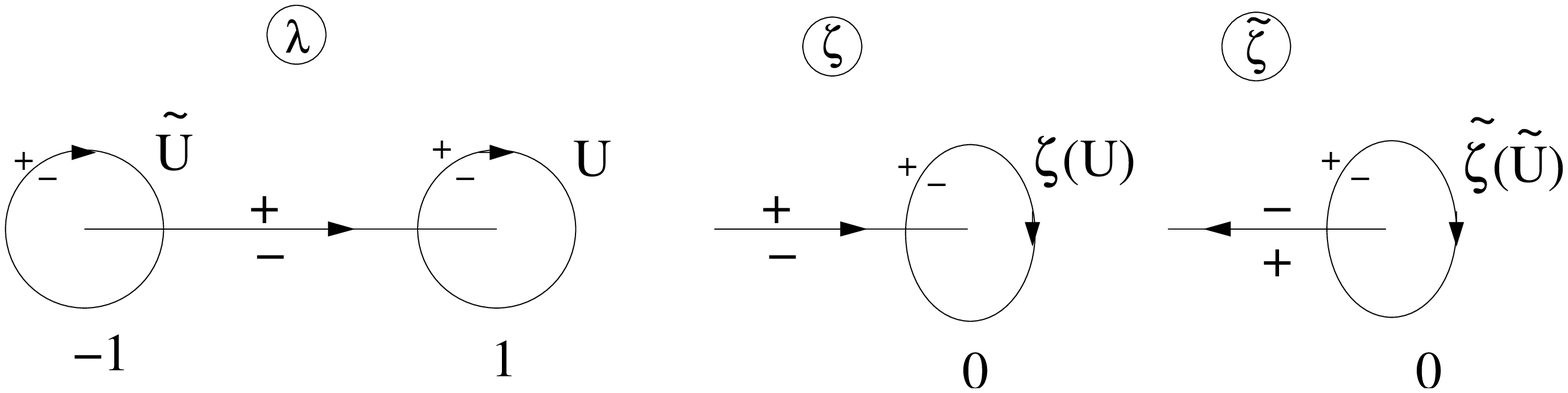,width=6.0in,angle=0}}
\vspace{0cm}
\caption{
Contours for the $Y$- and $\Phi$-RH problems and parametrices.}
\label{fig3}
\end{figure}

Consider the function
\be\la{f}
f(\lb)=\frac{1+i(\lb^2-1)^{1/2}\sin(\al/2)}{1-i(\lb^2-1)^{1/2}\sin(\al/2)}
\ee
which is analytic and has no zeros in $U\setminus (1-\de,1]$ and
$\wt U\setminus [-1,-1+\de)$. (Note, however, that it is singular
outside of these disks at $\lb=-i\cot(\al/2)$.) The branch of the root
is taken such that $(\lb^2-1)^{1/2}>0$ for $\lb>1$.
The function (\ref{f}) has the following boundary values on $(-1,1)$:
\be
f_+(x)=\frac{1-\sqrt{1-x^2}\sin(\al/2)}{1+\sqrt{1-x^2}\sin(\al/2)},
\qquad f_-(x)=f_+(x)^{-1},\qquad x\in(-1,1).\label{f+-}
\ee
Consider first the neighborhood $U$.
We look for a parametrix, an analytic function in
$U\setminus (1-\de,1]$, satisfying the jump condition of the
$\Phi$-RH problem
on $(1-\de,1)$, of the form
\be
P(\lb)=E(\lb)\hat P(\lb) f(\lb)^{-\si_3 n/2}, \qquad \lb\in
U\setminus (1-\de,1],
\ee
where $E(\lb)$ is a non-zero analytic matrix-valued function in $U$
(which therefore does not affect the jump condition)
to be chosen below so that $P$ matches $N$ to leading order on the boundary
$\partial U$.

It is easy to verify using (\ref{f+-}) that for $P$ to satisfy the jump
condition for the $\Phi$-RH problem across $(1-\de,1)$, $\hat P$ must
satisfy the jump relation
\be\la{jumphP}
\hat P_+(x)= \hat P_-(x)\pmatrix{0&1\cr -1&2}, \qquad x\in(1-\de,1).
\ee
An appropriate matrix function satisfying this jump relation was
constructed in \ci{DIZ} (cf. \ci{DIZ} (4.79), (4.871)).

For $\lb\in U\setminus (1-\de,1]$, define the analytic function
\be\la{om}
\om(\lb)={1\over 2}\ln f(\lb)
\ee
$$
\equiv i\sum_{k=0}^{\infty}(-1)^{k}\frac{1}{2k+1}
\sin^{2k+1}{\frac{\alpha}{2}}(\lb^2 -1)^{k+1/2}.
$$
Note that
\be\la{om1}
\om(\lb)=i\sqrt{2}\sin{\al\over
2}(\lb-1)^{1/2}G(\lb),
\ee
where $G(\lb)$ is analytic in all of
$U$, and
\be\la{G}
G(\lb)=1+(\lb-1)\left({1\over 4}-
{2\over 3}\sin^2{\al\over 2}\right)+O((\lb-1)^2),
\ee
for $\lb$ near $1$. Thus,
\begin{eqnarray}
e^{\om(\lb)}=f(\lb)^{1/2},\qquad \lb\in U\setminus (1-\de,1],\\
\om(x)_+=e^{i\pi}\om(x)_-,\qquad x\in(1-\de,1).
\end{eqnarray}
Furthermore, the function $\om^2(\lb)$ is analytic in all of $U$ and
\be\la{om_s}
\om^2(\lb)=2ue^{i\pi}\sin^2{\al\over 2}\left(1+{u\over 2}-
{4\over 3} u \sin^2{\al\over 2}+O(u^2)\right),\qquad u=\lb-1.
\ee
The term $O(u^2)$ in (\ref{G},\ref{om_s}) is uniform for all
$0\le\al\le\pi$. In fact, the estimate (\ref{om_s}) is uniform
for $\alpha$ belonging to any compact set in the complex
$\alpha$-plane. Let us choose $ 0< \ep < \pi$. Then,
for sufficiently small $\de$,
the asymptotic relation (\ref{om_s}) implies that 
\be\la{om3}
|\omega(\lambda)| \geq \sqrt{\de}\left|\sin{\al\over
2}\right|,\quad \forall \lb \in \partial U,\quad
\forall \al\in[0,\pi-\ep]\cup\cD_\ep(\pi).
\ee
Here $\cD_\ep(\pi)$ is the disk in the $\al$-plane of radius 
$\ep$ centered at $\al=\pi$. 

Introduce the new variable
\be
\ze=e^{-i\pi}n^2\om^2(\lb).\la{ze}
\ee
Note that the mapping $\lb\to\ze$ of $U$ is one-to-one.

>From (\ref{om_s}) and (\ref{om3}) it follows
that for $\de$ and $\ep$ sufficiently small, the following
inequalities hold:
\be\la{ineq1}
-\frac{3\pi}{4} \leq \arg{\sqrt{\ze}} \leq \frac{3\pi}{4},
\ee
and
\be\la{ineq2}
|\sqrt{\ze}| \geq n\sqrt{\de}\left|\sin{\al\over 2}\right| 
\equiv \rho\sqrt{\de},
\ee
$$
\forall \lb \in \partial U,\quad
\forall \al\in[0,\pi-\ep]\cup\cD_\ep(\pi).
$$
Inequality (\ref{ineq2}) together with (\ref{rhos05}) imply the estimate
\be\la{ze1}
|\sqrt{\ze}| \geq \frac{\sqrt{\de}}{2}s_0>1,
\ee
$$
\forall \lb \in \partial U,\quad
\forall \al\in\left[\frac{2s_0}{n},\pi-\ep\right]\cup\cD_\ep(\pi),\quad
\frac{2}{\sqrt{\de}} < s_0 \leq  n.
$$ 
A function $\hat P(\lb)$ analytic in $U\setminus (1-\de,1]$ and
satisfying (\ref{jumphP}) is given by the following expression in terms of
Hankel functions (cf. \ci{DIZ}) where $\sqrt{\ze}=e^{-i\pi/2}n\om(\lb)$:
\be
\hat P(\lb)=\pmatrix{H_0^{(1)}(\sqrt{\ze}) & H_0^{(2)}(\sqrt{\ze})\cr
\sqrt{\ze}\left(H_0^{(1)}\right)'(\sqrt{\ze}) &
\sqrt{\ze}\left(H_0^{(2)}\right)'(\sqrt{\ze})}.
\ee

Inequality (\ref{ineq1}) and estimate (\ref{ze1}) allow
us to use the standard expansion for Bessel functions and
obtain the following asymptotics on the boundary $\partial U$:
\be\eqalign{
\hat P(\lb)=
\sqrt{2\over\pi}\ze^{-\si_3/4}\pmatrix{1&1\cr i&-i}
\left( I+{i\over 8\sqrt{\ze}}\pmatrix{1&2\cr -2&-1}+
{3\over 2^7\ze}\pmatrix{1&-4\cr -4&1}+
\hat P_{r}(\lb)\right)\times\\
e^{n\om(\lb)\si_3}e^{-i(\pi/4)\si_3},
\qquad \lb \in \partial U,}\la{P_as}
\ee
where the remainder $\hat P_{r}(\lb)$ satisfies the uniform estimate
\be\la{P_as1}
|\hat P_{r}(\lb)| <\frac{C_{0}}{|\ze|^{3/2}},
\ee
$$
\forall \lb \in \partial U,\quad
\forall \al\in\left[\frac{2s_0}{n},\pi-\ep\right]\cup\cD_\ep(\pi),\quad
\frac{2}{\sqrt{\de}} < s_0 \leq n. 
$$
Here $C_{0}$ is a numerical positive constant which comes from
the universal asymptotic expansion of the Hankel function
$H_0^{(1)}(\sqrt{\ze})$ for $|\ze| > 1$ and 
$-3\pi/4 \leq \arg{\sqrt{\ze}} \leq 3\pi/4.$

Now let us choose $E(\lb)$ so that $P$ matches $N$ on $\partial U$
to leading order in $\rho$, i.e., $PN^{-1}\sim I$.
Clearly, we should take
\be
E(\lb)=N(\lb)e^{i(\pi/4)\si_3}{1\over 2}\pmatrix{1&-i\cr 1&i}
\sqrt{\pi\over 2}\ze^{\si_3/4}.
\ee
It is easy to verify that $E(\lb)$ has no branch point or singularity
at $\lb=1$. Hence $E(\lb)$ is analytic in $U$.

Thus, the parametrix in the neighborhood $U$ is given by the
expression:
\be\eqalign{
P(\lb)=
N(\lb)e^{i(\pi/4)\si_3}{\sqrt{\pi}\over 2^{3/2}}\pmatrix{1&-i\cr 1&i}
\ze^{\si_3/4}\times\\
\pmatrix{H_0^{(1)}(\sqrt{\ze}) & H_0^{(2)}(\sqrt{\ze})\cr
\sqrt{\ze}\left(H_0^{(1)}\right)'(\sqrt{\ze}) &
\sqrt{\ze}\left(H_0^{(2)}\right)'(\sqrt{\ze})}
f(\lb)^{-\si_3 n/2},}
\ee
where $\ze$ and $f(\lb)$ are defined by (\ref{ze},\ref{om},\ref{f}).

The construction of a parametrix in the neighborhood $\wt U$ is
similar. In this case, instead of (\ref{om}) we set
\be\la{omti}
\wt\om(\lb)=-{1\over 2}\ln f(\lb),
\ee
which is analytic in $\wt U\setminus [-1,-1+\de)$. Thus
\begin{eqnarray}
e^{-\wt\om(\lb)}=f(\lb)^{1/2},\qquad \lb\in \wt U\setminus [-1,-1+\de),\\
\wt\om(x)_+=e^{-i\pi}\wt\om(x)_-,\qquad x\in(-1,-1+\de).
\end{eqnarray}
We find the same power series expansion for $\wt\om^2(\lb)$ as
(\ref{om_s}) with $\lb-1$ replaced by $-\lb-1$.

We define the $\wt\ze$ variable for $\lb$ in $\wt U$ again by the equation
\be
\wt\ze=e^{-i\pi}n^2\wt\om^2(\lb).
\ee

Note that for both the images $\wt\ze(\wt U)$ and $\ze(U)$,
the slit for $\wt\ze$ (respectively $\ze$) variable lies
along the negative half-axis (if $\alpha$ is real;
it is slightly rotated away from the negative half-axis
if $\alpha$ is complex). However, the orientation is 
changed (see Figure 3).

With the above notation for $\wt\om$ and $\wt\ze$,
the parametrix in $\wt U$ matching $N(\lb)$ to leading order
at $\partial \wt U$ is given by the following expression:
\be\eqalign{
\wt P(\lb)=
N(\lb)e^{-i(\pi/4)\si_3}{\sqrt{\pi}\over 2^{3/2}}\pmatrix{i&1\cr -i&1}
\wt\ze^{-\si_3/4}\times\\
\si_1
\pmatrix{H_0^{(1)}(\sqrt{\wt\ze}) & H_0^{(2)}(\sqrt{\wt\ze})\cr
\sqrt{\wt\ze}\left(H_0^{(1)}\right)'(\sqrt{\wt\ze}) &
\sqrt{\wt\ze}\left(H_0^{(2)}\right)'(\sqrt{\wt\ze})}
\si_1
f(\lb)^{-\si_3 n/2},\qquad \si_1=\pmatrix{0&1\cr 1&0},}
\ee
where $\sqrt{\wt\ze}=e^{-i\pi/2}n\wt\om(\lb)$.

Following the steepest descent method, we now formulate a RH-problem for
the function
\be
R(\lb)=\cases{
\Phi(\lb)N(\lb)^{-1},& $\lb\in\complex\setminus(\overline{U\cup\wt U}
\cup(-1,1))$\cr
\Phi(\lb)P(\lb)^{-1},& $\lb\in U\setminus(1-\de,1]$\cr
\Phi(\lb)\wt P(\lb)^{-1},& $\lb\in\wt U\setminus [-1,-1+\de)$.}\la{Rdef}
\ee
By construction, the function $R(\lb)$ has no jumps across
$(1-\de,1)\cup(-1,-1+\de)$. Moreover, since {\it apriori} $R(\lambda)$
can have no stronger than logarithmic singularities at the
points $\pm 1$, the function $R(\lb)$ is in fact analytic in
the union of the discs $U\cup \wt U$. It solves the following   
RH-problem on the contour
$\Si=\partial U\cup\partial\wt U\cup(-1+\de,1-\de)$ (see Figure 4): 
\begin{itemize}
\item $R(\lambda)$ is holomorphic for all $\lambda \notin \Si$
\item $R(\infty) = I$
\item $R_{+}(\lambda) = R_{-}(\lambda)\Lambda(\lambda), \quad
\lambda \in \Si^{(0)} \equiv \Si \setminus \{1-\de, -1 + \de\}$, where
\beqa
\Lambda(x)=N_+(x)\pmatrix{1& -2f_+^n(x)\cr 0 &1}N_+(x)^{-1},
\quad x\in (-1+\de,1-\de),\la{RJ1}\\
\Lambda(\lambda)=P(\lb)N(\lb)^{-1},\quad \lb\in \partial U\setminus\{1-\de\},\la{RJ2}\\
\Lambda(\lambda)=\wt P(\lb)N(\lb)^{-1},\quad \lb\in \partial\wt U\setminus\{-1+\de\},
\la{RJ3}
\eeqa
\end{itemize}

\begin{figure}
\centerline{\psfig{file=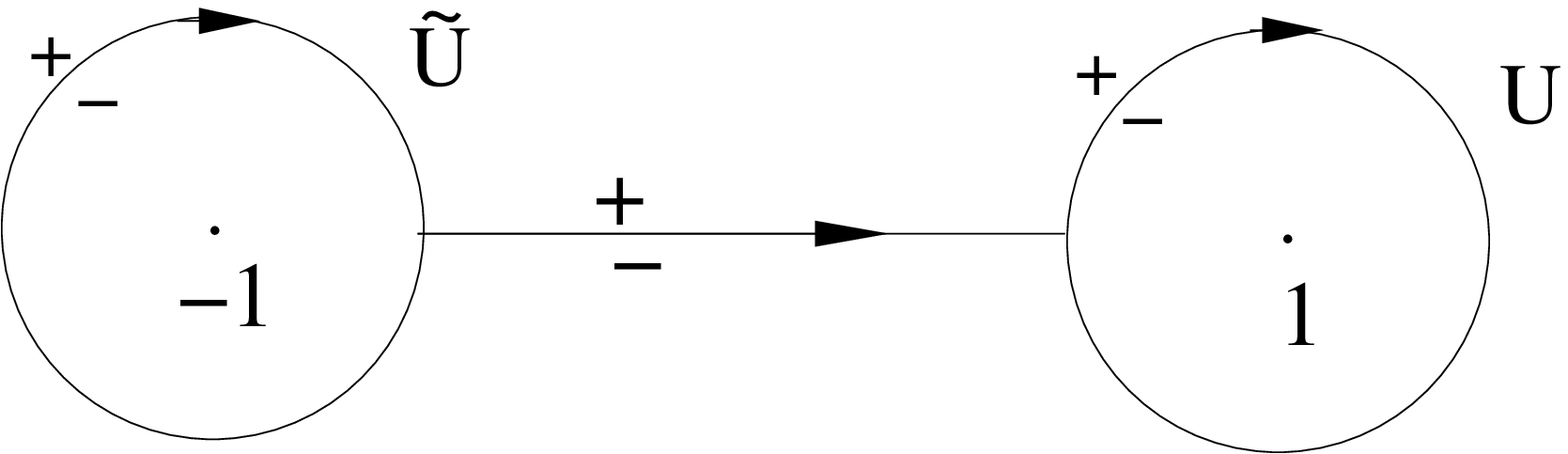,width=3.6in,angle=0}}
\vspace{0cm}
\caption{
Contour $\Si$ for the $R$-RH problem.}
\label{fig4}
\end{figure}

Observe that for all $-1 \leq x \leq 1$ and all $ 0 \leq \alpha \leq \pi$,
we have
$$
0 \leq f_{+}(x) \leq e^{-\sqrt{1-x^2}\sin{\al\over 2}}.
$$
Moreover, for sufficiently small $\ep$ there
exists a positive constant $C_{\de}$, depending on $\de$ only,
such that
$$
|f_{+}(x)| \leq e^{-C_{\de}},
$$
for all $-1 + \de \leq x \leq 1 - \de$ and all $\al \in \cD_\ep(\pi)$. 
Combining the two estimates above, we conclude that the
jump matrix on $[-1+\de,1-\de]$ is of order
$$
I+O(\exp{(-C_{\de, \ep}\rho)}),
$$ 
where $C_{\de, \ep}$ is a positive constant
which only depends on $\de$ and $\ep$.
This estimate is uniform in 
$$
x \in [-1+\de,1-\de],\quad 
\al \in [0,\pi-\ep]\cup\cD_\ep(\pi).
$$

Using (\ref{P_as}), we obtain the following asymptotic expansion
in inverse powers of $\sqrt{\ze}$ for the jump matrix
on $\partial U$:
\be\eqalign{
P(\lb)N(\lb)^{-1}=I+\La_1+\La_2+\Lambda_{r},\\
\La_1={i\over 16\sqrt{\ze}}\pmatrix{3\beta^2-\beta^{-2} &
  i(3\beta^2+\beta^{-2})\cr
i(3\beta^2+\beta^{-2}) & -(3\beta^2-\beta^{-2})},
\qquad
\La_2={3\over 2^7\ze}\pmatrix{1& -4i\cr 4i &1},\\
\lb\in\partial U,}\la{Lap1}
\ee
where $\beta(\lb)$ is defined in (\ref{Phiinfty}).
Since the matrix functions $N(\lambda)$ and $N^{-1}(\lambda)$
are uniformly bounded on $\partial U$, we conclude from 
(\ref{P_as}) and (\ref{P_as1}) that the error term $\Lambda_{r}(\lambda)$
in (\ref{Lap1}) satisfies the uniform estimate,
\be\la{P_as2}
|\Lambda_{r}(\lb)| <\frac{C_{\de}}{|\ze|^{3/2}},
\ee
$$
\forall \lb \in \partial U,\quad
\forall \al\in\left[\frac{2s_0}{n},\pi-\ep\right]\cup\cD_\ep(\pi),\quad
\frac{2}{\sqrt{\de}} < s_0 \leq n. 
$$
Here $C_{\de}$ is a positive constant depending on $\de$ only.
The jump matrix on $\partial\wt U$ is given by the similar representation with
the matrices $\La$ defined as follows:
\be\eqalign{
\La_1={i\over 16\sqrt{\wt\ze}}\pmatrix{-(3\beta^{-2}-\beta^{2}) &
  i(3\beta^{-2}+\beta^2)\cr
i(3\beta^{-2}+\beta^2) & 3\beta^{-2}-\beta^2},
\qquad
\La_2={3\over 2^7\wt\ze}\pmatrix{1& 4i\cr -4i &1},\\
\lb\in\partial\wt U.}\la{Lam1}
\ee

Let us summarize the above calculation.
{\proposition \la{prop1} The jump matrix $\Lambda(\lambda)$ 
of the $R$-RH problem
possesses the following properties:
\begin{enumerate}
\item For sufficiently small $\de$ and $\ep$,
      the function $\Lambda$ satisfies the estimates:
\begin{equation}\label{Las105}
|I - \Lambda(\lambda)| \leq \frac{C_{\de}}{\rho},
\quad \lambda \in (\partial U\setminus \{1-\de \})\cup 
(\partial \wt U\setminus \{-1+\de\}),
\end{equation}
and
\begin{equation}\label{Las205}
|I - \Lambda(x)| \leq \wt C_{\de, \ep}
\exp (- C_{\de, \ep}\rho),
\quad x \in (-1 +\de, 1 - \de),
\end{equation}
$$
\forall \al\in\left[\frac{2s_0}{n},\pi-\ep\right]\cup\cD_\ep(\pi),
\quad s_{0}\leq n.
$$
Here $\rho=n|\sin(\alpha/2)|$,
and $C_{\de}$, $\wt C_{\de, \ep}$, and $C_{\de, \ep}$
are positive constants depending on the indicated quantities only.
The number $s_{0}$ is any positive number satisfying the inequality
$s_{0}> 2/\sqrt{\de}$.
Moreover,
\begin{equation}\la{rho05}
\rho \geq \frac{s_0}{2},
\ee
$$
\forall \al\in\left[\frac{2s_0}{n},\pi-\ep\right]\cup\cD_\ep(\pi),\quad
s_{0} \leq  n.
$$
\item The estimate (\ref{Las105}) can be extended to the asymptotic series
\be \la{Las305}
\Lambda(\lambda) = I + \sum_{j=1}^{k-1}\Lambda_{j}(\lambda) + 
\Lambda^{(k)}_{r}(\lambda),
\quad \lambda \in  (\partial U\setminus \{1-\de \})\cup 
(\partial \wt U\setminus \{-1+\de\}),
\ee
where the terms $\Lambda_{j}$ of expansion (\ref{Las305}) and 
the error term $\Lambda^{(k)}_{r}(\lb)$ satisfy the uniform
estimates:
\be \la{Las405}
|\Lambda_{j}(\lambda)| \leq \frac{C^{(j)}_{\de}}{\rho^j}, \quad
|\Lambda^{(k)}_{r}(\lambda)| \leq \frac{C^{(k)}_{\de}}{\rho^k},
\ee
$$
\forall \al\in\left[\frac{2s_0}{n},\pi-\ep\right]\cup\cD_\ep(\pi),\quad
s_{0} \leq n. 
$$
The positive constants $C^{(j)}_{\de}$, $j =1,..., k$
are the functions of $\de$  only.
The first two terms of the expansion (\ref{Las305}), i.e. the
functions  $\Lambda_{1}$ and
$\Lambda_{2}$ are given by equations (\ref{Lap1}) if 
$\lambda \in  \partial U\setminus \{1-\de \}$, and 
 by equations (\ref{Lam1}) if 
$\lambda \in  \partial \wt U\setminus \{-1+\de \}$.
\item Let $\Lambda_{u}$,  $\Lambda_{d}$ , and $\Lambda_{l}$ denote the
limits of $\Lambda(\lambda)$ as $\lambda$ approaches the node point $1-\de$
from the above, from the below, and from the left along $\Si$, 
respectively. Then these limits exist, and the  following
cyclic equation holds:
\be\la{cycl105}
\Lambda_{d}\Lambda_{l}\Lambda^{-1}_{u} = I.
\ee
A similar relation (with $\Lambda_{l}$ replaced by
$\Lambda_{r}$) holds at the node point $-1+\de$.
\item The matrix function $\Lambda(\lambda)$ admits an analytic continuation
into a neighborhood of any point of the interval $(-1+\de, 1-\de)$. 
Moreover, this analytic continuation preserves the estimate (\ref{Las205})
with constants $\wt C_{\de,\ep}$, $C_{\de,\ep}$ possibly somewhat modified.
\end{enumerate} 
}
The only statements which need comments are the statements \# 3 and \#4.
These statements follow directly from the explicit formulae
(\ref{RJ1} - \ref{RJ3}) for the jump matrix $\Lambda(\lambda)$.

\bc \la{jamp05} The following inequalities hold:
\begin{eqnarray}
||I-\Lb||_{L^{2}(\Si)\cap L^{\infty}(\Si)}\le 
\frac{C^{(1)}_{\de, \ep}}{\rho},\la{jamp1051}\\
||I +\sum_{j=1}^{k-1}\Lambda_{j} - \Lambda ||_{L^{2}(\Si)\cap L^{\infty}(\Si)}
\leq  \frac{C^{(k)}_{\de, \ep}}{\rho^{k}},\la{jamp105}
\end{eqnarray}
$$
\forall \al\in\left[\frac{2s_0}{n},\pi-\ep\right]\cup\cD_\ep(\pi),\quad
s_{0} \leq  n, 
$$
where we set $\La_j\equiv 0$ for $\lb\in(-1+\de,1-\de)$.
\ec

By standard arguments of the $L^2$ RH theory (see e.g. \cite{DIZ,DZL2}), 
the inequality (\ref{jamp1051})
implies the solvability of the $R$-RH problem for sufficiently large $s_{0}$.
Moreover, let $\Omega_{k}$, $k =1, 2, 3$ denote the connected components
of the set ${\mathbb C}\setminus \Si$. Then, due to the cyclic relation 
(\ref{cycl105}), the restriction $R|_{\Omega_{k}}(\lambda)$ is continuous 
in $\overline{\Omega_{k}}$ for each $k$ (see, e.g., \ci{BDT}).

To complete the proof of the theorem, we need to show that the solution
$R(\lambda)$ of the $R$-RH problem satisfies the estimates indicated in
(\ref{th2105}).

{\lemma \la{lemma1} For sufficiently small $\de$ and $\ep$,
and for every $k$, the function $R(\lb)$ admits the asymptotic representation,
\be\label{Rexp1}
R(\lb)=I+\sum_{j=1}^{k-1}R_j(\lb)+R_r^{(k)}(\lb),
\ee
where
\be
R_j(\lb)=O\left( \frac{1}{(1+|\lambda|)\rho^j}\right),\qquad
R_r^{(k)}=O\left( \frac{1}{(1+|\lambda|)\rho^k}\right),\qquad
\rho \equiv n\left|\sin{\alpha \over 2}\right| \to \infty,
\ee
uniformly for all $\lb\in \Omega^{(2\de)}$
and $\al\in\left[\frac{2s_0}{n},\pi-\ep\right]\cup\cD_\ep(\pi)$. 
As in the Remark 1 to Theorem 2, the latter statement
means that there exist positive constants $C$ and $s_{0}$ such
that 
\be \la{newerr}
|R_{j}(\lambda)| \leq \frac{C}{(1+|\lambda|)\rho^j}, \quad
|R^{(k)}_{r}(\lambda)| \leq \frac{C}{(1+|\lambda|)\rho^k},
\ee
$$
\forall \lb\in\Omega^{(2\de)}, \quad
\forall \al\in\left[\frac{2s_0}{n},\pi-\ep\right]\cup\cD_\ep(\pi),\quad
\forall n: \quad s_{0} \leq n. 
$$
The functions $R_j(\lb)$
are constructed by induction as follows:
\beqa
R_1(\lb)={1\over 2\pi i}\int_{\partial U\cup\partial\wt U}
\Lb_1(s){ds\over s-\lb},\qquad
R_2(\lb)={1\over 2\pi i}\int_{\partial U\cup\partial\wt U}
(R_{1\,-}(s)\Lb_1(s)+\Lb_2(s)){ds\over s-\lb},\\
\dots,\qquad
R_{k-1}(\lb)={1\over 2\pi i}\int_{\partial U\cup\partial\wt U}
\sum_{j=1}^{k-1} R_{k-1-j,-}(s)\Lb_j(s){ds\over s-\lb},\qquad R_0\equiv I.
\eeqa
}

{\bf Remark 4} We also assume (cf. Remark 1 to Theorem 2) that $\ep$ is
small enough so that $\rho\ge s_0/2$
for all $\al\in\left[\frac{2s_0}{n},\pi-\ep\right]\cup\cD_\ep(\pi)$
and $n \geq s_{0}$.

The {\bf proof} of the lemma is essentially a combination of the arguments
from \ci{Dstrong} and \ci{KMVV}. We consider in detail
the case of $k=3$, which is all that is needed below, but the argument 
extends in an obvious way to any $k=1,2,\dots$. The details are left to 
the interested reader.

Write the jump condition for $R(\lb)$ in the form
\be\la{j}
R_{0\,+}+R_{1\,+}+R_{2\,+}+R_{r\,+}=(R_{0\,-}+R_{1\,-}+R_{2\,-}+R_{r\,-})
(I+\Lb_1+\Lb_2+\Lb_r).
\ee
Here $\Lb_1$, $\Lb_2$ are given by (\ref{Lap1},\ref{Lam1}) on 
$\partial U$, $\partial\wt U$,
and we set $\Lb_1=\Lb_2=0$ on $(-1+\de,1-\de)$. Thus $\Lb_r=O(1/\rho^3)$ on 
$\partial U\cup\partial\wt U$ (this error term arises from the Bessel 
asymptotics: see (\ref{Las405})),
and $\Lb_r=O(e^{-C_{\de,\ep}\rho})$ on  $(-1+\de,1-\de)$.
We now show that we can define $R_1$ and $R_2$ so that they are of order 
$1/\rho$ and $1/\rho^2$, 
respectively. We then show that the remainder $R_r$ is of order $1/\rho^3$.
Set 
\[
R_0=I.
\]
We define $R_j$ by collecting in 
(\ref{j}) the terms that we want to be of the same order. First,
\be\la{j1}
R_{1\,+}(\lb)=R_{1\,-}(\lb)+\Lb_1(\lb),\qquad \lb\in\Si.
\ee
We are looking for a function $R_1(\lb)$, which is holomorphic outside 
$\Si$, satisfying
$R_1(\lb)=O(1/\lb)$, $\lb\to\infty$, and the above jump condition.
The solution to this RH-problem is given by the Sokhotsky-Plemelj formula,
\be
R(\lb)=C(\Lb_1),
\ee
where
\[
C(f)={1\over 2\pi i}\int_{\Si}f(s){ds\over s-\lb}
\]
is the Cauchy operator on $\Si$.
The condition $\Lb_1(\lb)=O(1/\rho)$, $\lb\in\Si$, 
$\rho\to\infty$ (uniform in $\al$),
implies that there exist $c, d_{0}, s_0>0$ such that 
\be\la{R11}
|R_1(\lb)|\le \frac{c}{(1+|\lambda|)\rho}, \quad n \geq s_0,
\ee
uniformly
in $\al\in\left[\frac{2s_0}{n},\pi-\ep\right]\cup\cD_\ep(\pi)$ and 
$\lb$ satisfying ${\rm dist}(\lb,\Si)\ge d_{0}$.
Actually, this estimate is uniform for all $\lb\in\complex\setminus\Si$ 
up to $\Si$.
This can be shown either by direct calculation (see below) or by suitably 
deforming the contour $\Si$. Indeed, since
\be\la{Rminus105}
R_{1}(\lambda) = {1\over 2\pi i}\int_{\partial U\cup \partial \wt U}
\Lambda_{1}(s){ds\over s-\lb},
\ee
the estimate (\ref{R11}) holds for $\lambda$ 
up to the interval $(-1 + \de', 1-\de')$, for any $\de' > \de$.
Since $\Lambda_{1}(\lambda)$ is analytic in
the neighborhood of ${\partial U\cup \partial \wt U}$
(as, actually, are $\Lambda_{j}(\lambda)$ for all $j$),
the contour of integration in (\ref{Rminus105}) can be deformed
so that the estimate holds up to  ${\partial U\cup \partial \wt U}$
as well. It also should be observed that,
by the same deformation of the contour of integration 
in (\ref{Rminus105}), one obtains analytic continuations of
both the functions $R_{1+}(\lambda)$ and $R_{1-}(\lambda)$
in the neighborhood of the contour  ${\partial U\cup \partial \wt U}$. 
Moreover, this analytic
continuation preserves the estimate (\ref{R11}).

%The latter argument is the following.
%
%\begin{figure}
%\centerline{\psfig{file=sin5.eps,width=3.0in,angle=0}}
%\vspace{0cm}
%\caption{
%Deformed part of the contour $\Si$.}
%\label{fig5}
%\end{figure}
%
%
%
%Take $\lb\in\complex\setminus\Si$ such that 
%$\rm{dist}(\lb,\Si)<\ep$). Deform the contour as shown in Figure 5.
%Here $\wt\Si$ is the same as $\Si$ with the dotted part replaced by the
%semicircle of radius $\ep$. $\wt R(\lb)$ is defined as shown, and $J(\lb)$ is the
%analytic continuation of the
%jump matrix for $R$ on $\Si$. (It is easily seen that the continuation exists
%in a neighborhood of the original $\Si$.)
%$\wt R_1(\lb)$ satisfies the same Riemann-Hilbert problem as $R_1(\lb)$ but on
%the contour $\wt\Si$.
%The argument leading to (\ref{R11}) for $R_1$ and $\Si$ holds for $\wt R_1$ and
%$\wt\Si$ as well. Therefore (see the figure),
%\be\la{R1esti}
%|R_1(\lb)|=|\wt R_1(\lb)|\le c_1/\rho.
%\ee
%Analysis of the analytic continuation $J(\lb)$ shows that we can find the same
%$c_1$ for all $\al$, $\rho>\rho_0$, and all
%$\lb\in \complex\setminus\Si$ up to the boundary $\Si$.
%

Now define $R_2(\lb)$ by the jump condition
\be\la{j2}
R_{2\,+}(\lb)=R_{2\,-}(\lb)+R_{1\,-}(\lb)\Lb_1(\lb)+\Lb_2(\lb),\qquad \lb\in\Si,
\ee
together with the 
requirement of analyticity for $\lb\in\complex\setminus\Si$, and 
the condition $R_2(\lb)=o(1)$ for $\lb\to\infty$.
The solution to this RH-problem is 
\be\la{R21}
\eqalign{
R_2(\lb)=C(R_{1\,-}\Lb_1+\Lb_2)=\\
{1\over 2\pi i}\int_{\partial U\cup \partial \wt U}
(R_{1\,-}(s)\Lambda_{1}(s) + \Lambda_{2}(s)){ds\over s-\lb},
\qquad \lb \in \complex \setminus (\partial U\cup \partial \wt U).}
\ee
Using (\ref{R11}), the estimates $\Lb_j=O(1/\rho^{j})$, and the
analyticity of $R_{1-}$ and $\Lb_j$ in the neighborhood 
of  ${\partial U\cup \partial \wt U}$, we obtain in the same 
way as for $R_1$: for some $c>0$
\be\la{R22}
|R_2(\lb)|\le \frac{c}{(1+|\lambda|)\rho^2},\qquad 
\lb\in \complex\setminus\Si,\qquad n \geq s_{0}
\ee
with the same uniformity and analyticity properties in $\al$ and $\lb$.
Below in the proof, the same symbol $c$ will stand for various constants
independent of $\al$, $\lb$, and $n$.

Now from (\ref{j},\ref{j1},\ref{j2}) we obtain 
\be\la{jr}
R_{r\,+}(\lb)=M(\lb)+R_{r\,-}(\lb)\Lb(\lb),\qquad \lb\in\Si,
\ee
where 
\[
M\equiv R_{2\,-}\Lb_1+(R_{1\,-}+R_{2\,-})\Lb_2+(I+R_{1\,-}+R_{2\,-})\Lb_r.
\]

{\bf Remark} In the terminology of \ci{DZsobolev}, equation (\ref{jr}) is an 
inhomogeneous RH-problem of type 2. 

Since $R_r=R-I-R_1-R_2$, the matrix function $R_{r}(\lambda)$
is holomorphic outside $\Si$ and satisfies the condition
$R_r(\lb)=o(1)$ as $\lb\to\infty$. Therefore,
\be\la{Rr}
R_r(\lb)=C(M)+C(R_{r\,-}(\Lb-I)),\qquad \lb\in \complex\setminus\Si.
\ee
(It is worth mentioning that, by virtue of property \# 3 of the jump matrix
$\Lambda(\lambda)$ formulated in proposition \ref{prop1}, 
equation (\ref{Rr}) is consistent with the
absence of the singularities of the function $R_{r}(\lambda)$
at the node points $1-\de$ and $-1+\de$.) Equation (\ref{Rr}),
in turn, implies that
\be\la{Rr-}
R_{r\,-}(\lb)=C_-(M)+C_-(R_{r\,-}(\Lb - I)),\qquad \lb\in\Si,
\ee
where $C_-(f)=\lim_{\lb'\to \lb}C(f)$, as $\lb'$ approaches a point $\lb\in\Si$
from the ``$-$'' side of $\Si$. Now defining the operator
\[
C_{\Lb}(f)\equiv C_-(f(\Lb - I)),
\]
we represent (\ref{Rr-}) in the form
\be\la{RDe}
(I-C_{\Lb})(R_{r\,-})=C_-(M).
\ee
Because of the $L^{\infty}$ part of the estimate (\ref{jamp105}),
and the fact that 
$C_-$ is a bounded operator from $L^2(\Si)$ to $L^2(\Si)$,
it follows that
the operator norm $||C_\Lb||=O(1/\rho)$, and hence $I-C_{\Lb}$ is invertible by 
Neumann series for $s_{0}$ (and therefore $\rho$) sufficiently large. 
Thus (\ref{RDe}) gives
\be\la{Rrep05}
R_{r\,-}=(I-C_{\Lb})^{-1}(C_-(M)).
\ee\
Moreover, using the $L^2$ part of the estimate (\ref{jamp105}),
we conclude that $||C_{-}(M)||_{L^2(\Si)} = O(\rho^{-3})$. Together with (\ref{Rrep05}),
this yields the uniform estimate,
\be\la{RestL2}
 ||R_{r\,-}||_{L^2(\Si)} \leq \frac{c}{\rho^3},
\ee
$$
\forall \al\in\left[\frac{2s_0}{n},\pi-\ep\right]\cup\cD_\ep(\pi),
\quad n\geq s_{0}.
$$
Combining the estimate (\ref{RestL2}) with equation (\ref{Rr}), we can
complete the proof of the lemma as follows.

First, assuming that ${\rm dist}(\lb,\Si)\ge d_{0}$,
we immediately arrive at the estimate
\be\la{Rr105}
|C(M)(\lambda)| \leq \frac{c}{(1+|\lambda|)\rho^3}, \quad n \geq s_{0},
\ee
for the first term in the r.h.s. of (\ref{Rr}), and the estimate
$$
|C(R_{r\,-}(\Lambda - I))(\lambda)| \leq \frac{c}{1 + |\lambda|}
||R_{r\,-}||_{L^2(\Si)}||\Lambda - I||_{L^2(\Si)}\leq
$$
\be\la{Rr205}
\frac{c}{(1 + |\lambda|)\rho^4},
\quad n \geq s_{0},
\ee
for the second term. Both the estimates are uniform in
$\al\in\left[\frac{2s_0}{n},\pi-\ep\right]\cup\cD_\ep(\pi)$. Together they
yield the estimate
\be\la{Rr305}
|R_{r}(\lambda)| \leq \frac{c}{(1+|\lambda|)\rho^3}, \quad n\geq s_{0},
\ee
uniformly
in $\al\in\left[\frac{2s}{n},\pi-\ep\right]\cup\cD_\ep(\pi)$ and 
$\lb$ satisfying ${\rm dist}(\lb,\Si)\ge d_{0}$.

Second, we observe
that the matrix $\Lambda_{r}(\lambda)$ coincides with
the matrix $\Lambda(\lambda) - I$  on the interval  $(-1+\de, 1 -\de)$.
Hence, by property \# 4 of the matrix function $\Lambda(\lambda)$
(see proposition \ref{prop1}),
the matrix function $\Lambda_{r}(\lambda)$  admits an analytic
continuation in the neighborhood of
any point of the interval $(-1+\de, 1 -\de)$, and this continuation
preserves the estimate, $\Lambda_{r} = O(e^{-C_{\de,\ep}\rho})$. This means that,
by bending the segment $(-1+\de, 1 -\de)$ of the contour $\Si$ we can
extend $\lambda$ in the estimate (\ref{Rr105}) up to the interval
$(-1 + 2\de, 1-2\de)$. Using property \# 4 of the jump
matrix $\Lambda(\lambda)$ one more time, we can rewrite the second
term in the r.h.s. of  equation (\ref{Rr}) as
\be\la{Rr405}
\eqalign{
C(R_{r\,-}(\Lambda - I))(\lambda) =
 \frac{1}{2\pi i}\int_{\partial U \cup \partial \wt U}
R_{r\,-}(s)(\Lambda(s) - I)\frac{ds}{s-\lambda}+\\
\frac{1}{2\pi i}\int_{\gamma{(d)}}
R_{r}(s)(\Lambda(s) - I)\frac{ds}{s-\lambda},}
\ee
if $\lambda$ lies above the interval $(-1 + 2\de, 1-2\de)$, and as
\be\la{Rr505}
\eqalign{
C(R_{r\,-}(\Lambda - I))(\lambda) =
 \frac{1}{2\pi i}\int_{\partial U \cup \partial \wt U}
R_{r\,-}(s)(\Lambda(s) - I)\frac{ds}{s-\lambda}+\\
\frac{1}{2\pi i}\int_{\gamma^{(u)}}
(R_{r}(s) - M(s))(I - \Lambda^{-1}(s))\frac{ds}{s-\lambda},}
\ee
if $\lambda$ lies below the interval $(-1 + 2\de, 1-2\de)$.
Here, the contours $\gamma^{(d)}$ and $\gamma^{(u)}$ are 
the slight deformations of the segment $(-1 + \de, 1-\de)$
down and up, respectively. Using, 
in representations (\ref{Rr405})
and (\ref{Rr505}), the estimate (\ref{Rr305}) for  $R_{r}(\lambda)$,
we extend the variable $\lambda$ 
in the  estimate (\ref{Rr205}) up to
the interval $(-1 + 2\de, 1-2\de)$. 

The above extensions of the estimates (\ref{Rr105}) and (\ref{Rr205})
mean, in particular, that they both, and hence the 
estimate (\ref{Rr305}), are valid for all $\lambda \in \Omega^{(2\de)}$.
The proof of the lemma is completed. $\Box$
\vskip .2in

We now derive explicit formulae for
the terms $R_{1}(\lambda)$ and $R_{2}(\lambda)$ of the expansion
(\ref{Rexp1}). By Lemma 1,
\be
R_1(\lb)={1\over 2\pi i}\int_{\partial U\cup\partial\wt U}
{\La_1(x)dx\over x-\lb},\qquad
R_2(\lb)={1\over 2\pi i}\int_{\partial U\cup\partial\wt U}
{R_{1\,-}(x)\La_1(x)+\La_2(x)\over x-\lb}dx\la{plem}
\ee

As noted in \ci{KMVV},
we can also obtain the expressions for $R_j(\lb)$
in the following way. It is not difficult to check that
$\La_1(\lb)$ and $\La_2(\lb)$ are analytic in $ (U\setminus \{1\}) \cup 
(\wt U\setminus \{-1\})$
with the simple poles at $\pm1$. We have
\be
\La_1(\lb)=\frac{A^{(1)}}{\lb-1}+O(1),\quad\mbox{as } \lb\to 1,
\qquad
    \La_1(\lb)=\frac{B^{(1)}}{\lb+1}+O(1),\quad\mbox{as } \lb\to -1,
\ee
where the constant matrices $A^{(1)}$ and $B^{(1)}$ are obtained by expanding
$\om(\lb)$ and $\beta(\lb)$ in (\ref{Lap1},\ref{Lam1}) at $\lb=\pm1$.
It is easy to verify directly that
the Riemann-Hilbert problem for $R_1(\lb)$ has the solution:
\be
R_1(\lb)=
\cases{\frac{A^{(1)}}{\lb-1}+\frac{B^{(1)}}{\lb+1},&
for $\lb\in\complex\setminus(\overline{U\cup\wt U})$\cr
\frac{A^{(1)}}{\lb-1}+\frac{B^{(1)}}{\lb+1}-\La_1(\lb),&
for $\lb\in U\cup\wt U$.}\la{R1}
\ee

Using the series (\ref{om_s}) and the expansion of $\beta(\lb)$ at $\pm1$,
it is not difficult to obtain the singular and constant term in the
Laurent expansion of $\La_1(\lb)$. By the first formula in (\ref{R1}),
we obtain (using the singular term) the expression (\ref{R105}).

Similarly we may calculate the singular term in the expansion of $\La_2(\lb)$ at
$\pm1$, and use the second formula in (\ref{R1}) to evaluate
$R_{1}(\pm1)$ (note that the formula (\ref{R105})
is valid only outside $U\cup\wt U$). 
It is then easy to compute the integral for $R_2$ in
(\ref{plem}) and obtain (\ref{R205}).
This completes the proof of the theorem.$\Box$.
\vskip .2in

Now we give some remarks and corollaries of Theorem \ref{theorem2}.

\noindent
{\bf Remark 5} Estimate (\ref{th21}) and formula (\ref{Phiinfty})
imply that
$$
\Theta(n,\alpha) \sim \cos{\frac{\alpha}{2}},
\quad \Delta(n,\alpha) \sim \sin^2{\frac{\alpha}{2}},
$$
and using either (\ref{differencePhi}) or (\ref{differentialPhi})
we recover the master term of Widom's asymptotics \ci{W1} (cf. also
\ci{DIZ}),
$$
\ln D_{n}(\alpha) \sim n^2\ln \cos{\frac{\alpha}{2}},
\quad n\to \infty.
$$
\vskip .2in

{\corollary \la{cor1} The function $\Delta(n,\alpha)$ admits
the asymptotic expansion 
\be\label{Delta_as}
\Delta(n,\al)=\sin^2{\al\over 2}-{\cos^2(\al/2)\over 4n^2}
+O\left({1\over\rho^3}\right)\sin^2{\al},\qquad \rho\to\infty,
\ee
which is uniform for $\al\in\left[\frac{2s_0}{n},\pi\right]$.
}

\noindent
{\bf Remark 6} The statement,
$$
\Delta_{r}(n,\alpha) = 
O\left({1\over\rho^3}\right)\sin^2{\al},\qquad \rho\to\infty
$$
uniformly for 
$ \al\in\left[\frac{2s_0}{n},\pi \right],$
means that there exist positive constants $C$ and $s_{0}$,
such that
\begin{equation}\label{meaningD05}
|\Delta_{r}(n,\alpha)| \leq \frac{C}{\rho^{3}}\sin^2{\alpha},
\end{equation}
$$
\forall \al\in\left[\frac{2s_0}{n},\pi\right],
\quad \mbox{and}\quad
s_{0} \leq  n.
$$

\noindent
{\bf Proof of Corollary 2.}
To calculate $\De$ we need the asymptotics of $\Phi(\lb)$ outside
the neighborhoods $U$ and $\wt U$. By (\ref{th2105}) these are given
by the expression:
\be\la{Phi}
\Phi(\lb)=(I+R_1+R_2+R_r^{(3)})N,
\qquad \lb\in\Omega^{(2\de)}
\ee
where $R_r^{(3)}$ is estimated by (\ref{newerr}) for $k=3$.
In particular, the estimate (\ref{newerr})  becomes
\be
O(\rho^{-3})\sin{\al\over2}
\qquad {\mathrm{if}}\quad \lb=\pm i\cot{\al\over 2}.
\la{eest}
\ee
Similarly,
\be
R_j(\pm i\cot{\al/2})=O(\rho^{-j})\sin{\al\over 2},
\qquad j=1,2,\dots \la{Rest}
\ee
Since
\be
N(\pm i\cot{\al/2})=\pmatrix{\cos(\al/4)&\pm\sin(\al/4)\cr
\mp\sin(\al/4)&\cos(\al/4)},
\ee
we have $N(-i\cot{\al/2})^{-1} = N(i\cot{\al/2}) = O(1)$ and
$[N(-i\cot{\al/2})^{-1}N(i\cot{\al/2})]_{12}=\sin{\al/2}$.
Definition (\ref{Deltadef}) and  equations (\ref{Phi},\ref{eest},\ref{Rest})
then imply
\begin{eqnarray}
\eqalign{
\Delta(n, \alpha) =
\left[N(-i\cot{\al/2})^{-1}\left(I+(O(\rho^{-1})+O(\rho^{-2})+O(\rho^{-3}))
\sin{\al\over2}\right)\times\right.\\
\left.
\left(I+(O(\rho^{-1})+O(\rho^{-2})+O(\rho^{-3}))
\sin{\al\over2}\right)N(i\cot{\al/2})\right]_{12}^2=}\label{deltamed}\\
=\Delta_{0}(\alpha) + \frac{1}{n}\Delta_{1}(\alpha)
+ \frac{1}{n^2}\Delta_{2}(\alpha) + {f(\al,n)\over \rho^3}\sin^2(\al/2),\label{deltaest}
\end{eqnarray}
where
\begin{equation}\label{Delta0}
\Delta_{0}(\alpha) = \sin^{2}{\frac{\alpha}{2}},
\end{equation}
and $f(\al,n)$ is uniformly bounded for 
$\al\in[2s_0/n,\pi-\ep]\cup\cD_\ep(\pi)$, and $s_{0} \leq n$.
Note that to write (\ref{deltamed}) we used the fact that $\det R(\lb)=1$.

In order to determine the terms $\Delta_{1}(\alpha)$
and $\Delta_{2}(\alpha)$ in this equation we
need $R_{1,2}(\pm i\cot{\al/2})$. These values we obtain from 
(\ref{R105},\ref{R205}):
\begin{eqnarray}
R_{1}(\pm i\cot{\al/2})=
\pm{1\over 8n}\pmatrix{-\cos(\al/2)&\pm\sin(\al/2)\cr
\pm\sin(\al/2)&\cos(\al/2)},\\
R_{2}(\pm i\cot{\al/2})=
\pm{1\over 2^7n\rho}\pmatrix{\pm\sin(\al/2)&-8\cos(\al/2)\cr
8\cos(\al/2)&\pm\sin(\al/2)}.
\end{eqnarray}
As $\det\Phi(\lb)=1$, the inverse $\Phi^{-1}$ is easy to compute, and
after a simple computation we arrive at the equations,
\be\label{Delta1_as}
\Delta_{1}(\al)= 0, \quad \Delta_{2}(\al) = -\frac{1}{4}\cos^2(\al/2).
\ee

Now it only remains to show that the function $f(\al, n)$
in (\ref{deltaest}) satisfies the estimate
\be\la{fest05}
|f(\al, n)| \leq C\cos^{2}(\al/2),
\ee
for $\al\in[2s_0/n,\pi]$, and $s_{0} \leq n$.
In fact, since the uniform boundedness of $f(\al, n)$
on the sets indicated
has already been established, it is enough to show that estimate (\ref{fest05})
holds for all $\al \in \cD_{\ep/4}(\pi)$.

Observe that, for fixed $n$, the quantity $\Delta(n,\alpha)$ is an
analytic function at $\alpha = \pi$ with $\Delta(n,\pi)=1$ and
$(d/d\al)\Delta(n,\pi)=0$ (see (\ref{Dn1}) and (\ref{differentialPhi}))
so that we can write down the Taylor series for $\De$ in $\beta = \pi -\alpha$:
\begin{equation}\label{deltataylor}
\Delta(\alpha)=1 + (\al -\pi)^2 a_{2} +\cdots.
\end{equation}
As follows from  equation (\ref{deltaest}),
$f(\al,n)$ is a holomorphic function of $\al$ in 
$\cD_\ep(\pi)$.
Using a representation of $f(\al,n)$ by a Cauchy integral, we obtain:
\be\eqalign{
f(\al,n)={1\over 2\pi i}\int_{\partial\cD_{\ep/2}(\pi)}
{f(\wt\al,n)\over \wt\al-\al}d\wt\al=
{1\over 2\pi i}\int_{\partial\cD_{\ep/2}(\pi)}
{f(\wt\al,n)\over \wt\al-\pi}d\wt\al+\\
{(\al-\pi)\over 2\pi i}\int_{\partial\cD_{\ep/2}(\pi)}
{f(\wt\al,n)\over(\wt\al-\pi)^2}d\wt\al+
{(\al-\pi)^2\over 2\pi i}\int_{\partial\cD_{\ep/2}(\pi)}
{f(\wt\al,n)\over(\wt\al-\pi)^2(\wt\al-\al)}d\wt\al,\qquad |\pi-\al|<\ep/4.}
\la{cauchy}
\ee
At the same time, from (\ref{deltaest}), (\ref{deltataylor}), and 
(\ref{Delta1_as}) it follows that the Taylor series of $f(\al,n)$
at $\al = \pi$ has the form,
$$
f(\al, n) = (\al - \pi)^2\tilde{a}_{2}+\dots
$$
Therefore, the first two integrals in the r.h.s. of (\ref{cauchy}) must be zero, 
and the third one, by virtue of the uniform
boundedness of $f(\al,n)$ for all 
$\al \in \cD_\ep(\pi)$ and all $n \geq s_{0}$,  yields the estimate
(\ref{fest05}) for all $\al \in \cD_{\ep/4}(\pi)$ and all
$n \geq s_{0}$. The proof of the corollary is completed. $\Box$

\vskip .2in
\noindent
{\bf Remark 7} Here is an alternative derivation of the leading terms
in formula (\ref{Delta_as}).
\vskip .2in

We start with equations (\ref{deltaest},\ref{Delta0}).
The issue is the exact evaluation of the quantities
$\Delta_{1}(\alpha)$ and $\Delta_{2}(\alpha)$. This can be done with
the help of the relevant (integrable) differential system associated in the
standard way with the original $m$ - RH problem. Indeed, it is shown in
\ci{DIZ} that the Toeplitz determinant $D_{n}(\alpha)$, considered as the
function of the variable
$$
t = e^{-2i\alpha}
$$
is the $\tau$ - function  for the  Painlev\'e VI equation
characterized by the parameters
$$
\theta_{\infty} = - \theta_{0} = n, \quad
\theta_{1} = \theta_{t} = 0,
$$
where we use the $\theta$-notations of Jimbo,
see \ci{Jimbo82}. According to \ci{Jimbo82}, this means that the quantity
\begin{equation}\label{etadef}
\eta(t) \equiv t(t-1)\frac{d}{dt}\ln D_{n}
\end{equation}
satisfies the following nonlinear differential equation (the $\tau$-
form of Painlev\'e VI):
$$
\left(\frac{d\eta}{dt} - \frac{n^2}{4}\right)
\left(t(t-1)\frac{d^{2}\eta}{dt^{2}}\right)^{2}
$$
\begin{equation}\label{painleve}
+ \left[2\left(\frac{d\eta}{dt} - \frac{n^2}{4}\right)
\left(t\frac{d\eta}{dt} - \eta\right)
-\left(\frac{d\eta}{dt}\right)^{2} +\frac{n^2}{2}\frac{d\eta}{dt}\right]^2
= \left(\frac{d\eta}{dt}\right)^{4}.
\end{equation}
The functions $\Delta(n, \alpha)$  and $\eta(t)\equiv \eta(n, t)$ are related
by the equation
\begin{equation}\label{Deltaeta}
\Delta = \frac{1-t}{n^2}\, \,\frac{d\eta}{dt} + \frac{1}{n^2}\eta,
\end{equation}
and we may anticipate an expansion for $\eta$ similar to (\ref{Delta_as}).
Indeed we expect
\begin{equation}\label{etaest}
\eta(t)\equiv \eta(n, t) = n^2\,\eta_{0}(t) + n\,\eta_{1}(t)
+ \eta_{2}(t) + O\left(\frac{1}{\rho}\right),
 \quad \rho \to \infty,
\end{equation}
where
\begin{equation}\label{eta0}
\eta_{0}(t) = \frac{1}{4}(1-\sqrt{t})^2.
\end{equation}
A substitution of the asymptotics (\ref{etaest}) into the
equation (\ref{painleve}) gives us, after a straightforward
calculation, the following formulae for the coefficient
functions $\eta_{1}(t)$ and $\eta_{2}(t)$:
\begin{equation}\label{eta12}
\eta_{1}(t) \equiv 0,\quad \eta_{2}(t) =
 - \frac{1}{16}(1+\sqrt{t})^2.
\end{equation}
These equations together with (\ref{Deltaeta})
lead immediately to the leading terms in the formula (\ref{Delta_as}).

It also should be noticed that the differentiability of 
the asymptotics (\ref{deltaest}) follow from its uniformity in
the disk $\cD_\ep(\pi)$.

\section{Asymptotic evaluation of $D_{n}(\alpha)$. Proof of
estimate (\ref{asnodelta}).}
The asymptotic evaluation of the Toeplitz determinant $D_{n}(\alpha)$
is based on the integration of the differential identity
(\ref{differentialPhi}) from $\al$ to $\al_0$ (which is
close to $\pi$ from below).
We have:
\be\la{Dinteg}
(\al_0-\al)(\ln D_n)'(\al_0)- \ln D_n(\al_0) +\ln D_n(\al)=
-n^2\int_\al^{\al_0}d\theta \int_\theta^{\al_0} {\De(\phi)\over \sin^2\phi}
d\phi.
\ee

Fix $n$ and set $\al_0=\pi-\beta$. Substituting for $\ln D_n(\pi-\beta)$
the expansion (\ref{Dn1}), and for $\De(\phi)$ the asymptotics
(\ref{Delta_as}), and after taking the limit $\beta\to 0$, we
immediately obtain (\ref{asnodelta}) with the remainder
$O(1/\{n\sin(\al/2)\})$ uniformly for
$\frac{2s_0}{n} \le \alpha \le \pi-\ep$, $n \geq s_{0}$, $\ep>0$.

\vskip .2in
 
{\bf Acknowledgements.}
\vskip .2in
Percy Deift was supported
in part by NSF grants \# DMS-0296084 and \# DMS 0500923. 
Alexander Its was supported
in part by NSF grants \# DMS-0099812 and \# DMS-0401009. Xin Zhou was supported
in part by NSF grant \# DMS-0071398.


\begin{thebibliography}{99}

\bibitem{M} M. L. Mehta: Random matrices. San Diego: Academic
  1990

\bibitem{dCM} J. des Cloizeaux and M. L. Mehta, 
Asymptotic behavior of spacing distributions for the eigenvalues of
random matrices, {\it J. Math. Phys.} {\bf14}, 1648--1650 (1973)

\bibitem{D} F. Dyson, Fredholm determinants and inverse
  scattering problems. {\it Commun. Math. Phys.} {\bf 47}, 171--183 (1976)

\bibitem{W1} H. Widom, The strong Szeg\H o limit theorem for
  circular arcs.  {\it Indiana Univ. Math. J.} {\bf 21}, 277--283 (1971)

\bi{W2} H. Widom, The asymptotics of a continuous analogue
  of orthogonal polynomials. {\it J. Approx. Th.} {\bf 77}, 51--64 (1994)

\bi{W3} H. Widom,
Asymptotics for the Fredholm determinant of the sine kernel 
on a union of intervals, {\it Comm. Math. Phys.} {\bf 171},  159--180 (1995) 

\bi{DIZ} P. Deift, A. Its, and X. Zhou, A Riemann-Hilbert
  approach to asymptotic problems arising in the theory of random
  matrix models, and also in the theory of integrable statistical
  mechanics. {\it Ann. Math} {\bf 146}, 149--235 (1997)

\bi{K} I. V. Krasovsky, Gap probability in the spectrum of random matrices
and asymptotics of polynomials orthogonal on an arc of the unit
circle. {\it Int. Math. Res. Not.} {\bf 2004}, 1249--1272 (2004)

\bi{E} T. Ehrhardt, Dyson's constant in the asymptotics of the Fredholm 
determinant of the sine kernel. arXiv.org:math.FA/0401205

\bi{BE} E. L. Basor, T. Ehrhardt, On the asymptotics of certain 
Wiener-Hopf-plus-Hankel determinants. arXiv.org:math.FA/0502039

\bi{BT} E. L. Basor, C. A. Tracy, Some problems associated with the asymptotics
of $\tau$-functions. Surikagaku (Mathematical Sciences) {\bf 30}, no. 3, 71--76
 (1992) [English translation appears in RIMS-845 preprint]

\bi{BB} A. M. Budylin, V. S. Buslaev, Quasiclassical asymptotics of the resolvent of 
an integral convolution operator with a sine kernel on a finite interval. 
(Russian)  Algebra i Analiz  {\bf 7},  no. 6, 79--103 (1995);  
translation in  St. Petersburg Math. J.  {\bf 7},  no. 6, 925--942 (1996) 


\bi{bdj} J. Baik, P. Deift, K. Johansson, On the distribution of the
length of the longest increasing subsequence of random
permutations, {\it J.\ Amer.\ Math.\ Soc.} {\bf 12}, No.~4, 1119--1178 (1999) 


\bi{deift} P. A. Deift, Integrable systems and combinatorial theory,
{\it Notices of the Ammer. Math. Soc.}, {\bf 47} (6), 631--640 (2000)

\bi{DIK2} P. Deift, A. Its, I. Krasovsky, Asymptotics of 
the Airy-kernel determinant. In preparation.

\bi{Deift} P. Deift, Orthogonal polynomials and random matrices: a
Riemann-Hilbert approach. Courant Lecture Notes in Math. 1998

\bibitem{Szego} G. Szeg\H o, Orthogonal polynomials. AMS Colloquium
  Publ. {\bf 23}. New York: AMS 1959

\bi{DZ} P. Deift and X. Zhou, A steepest descent method for
oscillatory Riemann-Hilbert problem. {\it Ann. Math.} {\bf 137}, 295--368
(1993)

\bi{DZpainleve2} P. Deift and X. Zhou,
Asymptotics for the Painlevé II equation. {\it Comm. Pure Appl. Math.}  {\bf 48},
277--337  (1995) 

\bi{DVZzerodisp} P. Deift, S. Venakides, and X. Zhou, New results in small 
dispersion KdV by an extension of 
the steepest descent method for Riemann-Hilbert problems.  
{\it Int. Math. Res. Not.} {\bf 1997}, 286--299 (1997). 

\bi{Dstrong} P. Deift, T. Kriecherbauer, K. T-R McLaughlin,
S. Venakides, X. Zhou, Strong asymptotics for orthogonal polynomials
with respect to exponential weights. {\it Commun. Pure Appl.Math.} {\bf 52},
1491--1552 (1999)

\bibitem{DZL2}
P. Deift and X. Zhou, A priori $L^{p}$ estimates for solutions
of Riemann-Hilbert problems, {\it Int.\ Math.\ Res.\ Noties},
{\bf 40}, 2121--2154 (2002).

\bibitem{BDT} R. Beals, P. Deift, C. Tomei, Direct and inverse scattering on 
the line. Mathematical Surveys and Monographs, {\bf 28}. AMS, Providence, RI, 1988.

\bibitem{DZsobolev} P. Deift and X. Zhou, Long-time asymptotics for solutions 
of the NLS equation with initial data in a weighted Sobolev space. 
{\it Comm. Pure Appl. Math.} {\bf 56}, 1029--1077 (2003).

\bi{LS} G. S. Litvinchuk and I. M. Spitkovskii, Factorization of measurable matrix functions,
Birkh\"auser, 1987

\bi{Zhou} X. Zhou, The Riemann-Hilbert problem and inverse scattering,
{\it SIAM J. Math. Anal.}  {\bf 20}, No.~4, 966--986 (1989) 


\bi{KMVV} A. B. J. Kuijlaars, K. T-R McLaughlin, W. Van Assche,
M. Vanlessen, The Riemann-Hilbert approach to strong asymptotics for
orthogonal polynomials on $[-1,1]$. {\it Adv. Math.} {\bf 188}, 337--398 (2004) 

\bi{Jimbo82} M. Jimbo, Monodromy problem and the boundary condition for
some Painleve equations, Publ. RIMS, Kyoto University 18, 1137-1161 (1982)

\end{thebibliography}
\end{document}